\documentclass[article,11pt]{jsg}

\usepackage{amsmath}
\usepackage{amssymb} 
\usepackage{amsfonts}
\usepackage{latexsym} 
\usepackage{epic,eepic,epsfig,psfrag} 
\usepackage{amscd}  

\usepackage{lscape}

\def\one{\hbox{1\hskip-2.7pt l}}
\def\smone{{\scriptstyle\rm 1\hskip-2.05pt l}}

\def\p{\phi}

\def\a{\alpha}
\def\b{\beta} 
\def\d{\delta} 

\def\ep{\varepsilon} 
\def\e{\eta} 

\def\th{\theta}

\def\t{\tau} 
\def\i{\iota} 
\def\x{\xi} 
 
\def\n{\nu} 
\def\m{\mu} 
 
\def\o{\omega} 
\def\l{\lambda} 
 
\def\L{\Lambda}

\def\S{\Sigma}
\def\Si{\Sigma} 
 
\def\X{\Xi} 
\def\Om{\Omega} 
\def\P{\Phi} 
\def\cA{{\mathcal A}}
\def\cB{{\mathcal B}} 
\def\cC{{\mathcal C}} 
 
\def\cE{{\mathcal E}} 
\def\cF{{\mathcal F}} 
\def\cG{{\mathcal G}}

\def\cL{{\mathcal L}} 
\def\cM{{\mathcal M}}

\def\cR{{\mathcal R}}

\def\cU{{\mathcal U}}

\def\cg{{\mathfrak g}}

\def\rd{{\rm d}}
\def\rT{{\rm T}}

\def\rG{{\rm G}}
\def\R{{\mathbb R}}

\def\N{{\mathbb N}} 
 
\def\H{{\mathbb H}} 
\def\Z{{\mathbb Z}}

\def\tu{{\tilde u}}

\def\tA{{\tilde A}}
\def\tB{{\tilde B}}

\def\tX{{\tilde\Xi}}
\def\dr{{\rm d}r}
\def\ds{{\rm d}s}
\def\dt{{\rm d}t}

\def\dph{{\rm d}\p}

\def\st{\: \big| \:} 
\def\la{\langle\,}
\def\ra{\,\rangle}
\def\comp{\circ}
\def\half{{\textstyle{\frac 12}}} 

\def\tint{\textstyle\int}
\def\pd{\partial}
\def\laplace{\Delta}

\def\im{{\rm im}\,}
\def\hol{{\rm hol}}
\def\tr{{\rm tr}}
\def\Hom{{\rm Hom}}

\def\HF{{\rm HF}}
\def\SO{{\rm SO}}
\def\SU{{\rm SU}}
\def\su{{\mathfrak s\mathfrak u}}
\def\CS{{\mathcal C}{\mathcal S}}
\def\Af{{\cA_{\rm flat}}}

\newtheorem{dfn}{Definition}[section] 
\newtheorem{lem}[dfn]{Lemma} 
\newtheorem{prp}[dfn]{Proposition} 
\newtheorem{thm}[dfn]{Theorem} 
 
\newtheorem{cor}[dfn]{Corollary} 
 
\newtheorem{con}[dfn]{Conjecture} 
        
\renewcommand{\labelenumi}{(\roman{enumi})}
\renewcommand{\labelenumii}{\alph{enumii})}

\begin{document}

\bibliographystyle{plain}

\author{Katrin Wehrheim \\ IAS Princeton}

\title[Lagrangian boundary conditions for ASD instantons]
{Lagrangian boundary conditions for anti-self-dual instantons 
 and the Atiyah-Floer conjecture}

\maketitle

\section{Introduction}

The purpose of this survey is to explain an approach to the Atiyah-Floer conjecture
via a new instanton Floer homology with Lagrangian boundary conditions.
This is a joint project with Dietmar Salamon; see \cite{Sa} for an earlier exposition.
This paper also provides a rough guide to the analysis of
anti-self-dual instantons with Lagrangian boundary conditions in 
\cite{W elliptic, W bubbling}, which is the crucial ingredient of our approach.

Atiyah \cite{A} and Floer conjectured a natural isomorphism between the instanton
Floer homology $\HF_*^{\rm inst}(Y)$ of a homology $3$-sphere $Y$ 
and the symplectic Floer homology $\HF_*^{\rm symp}(\cR_\S,L_{H_0},L_{H_1})$ 
of a pair of Lagrangians $L_{H_0},L_{H_1}$ in the symplectic moduli space 
$\cR_\S$ of flat $\SU(2)$-connections, associated to a Heegard splitting 
$Y=H_0\cup_\S H_1$.
Both homologies were introduced by Floer \cite{F1,F2}, but the symplectic 
Floer homology is not strictly defined in this case due to singularities of $\cR_\S$.
Taubes \cite{T} proved that the Euler characteristics both agree
with the Casson invariant of $Y$.
The main task in identifying the homology groups is a comparison between the
trajectories:
pseudoholomorphic curves in $\cR_\S$ with Lagrangian boundary conditions 
and anti-self-dual instantons on $\R\times Y$ (which has no boundary).

The basic idea of our approach is to introduce a third 
Floer homology\footnote{This is a 
special case of the invariant $\HF_*^{\rm inst}(Y,\cL)$ introduced below
for a $3$-manifold $Y$ with boundary and a Lagrangian submanifold $\cL$
in the space of connections over $\pd Y$.
}${\HF_*^{\rm inst}([0,1]\times\S,\cL_{H_0}\times\cL_{H_1})}$
whose trajectory equation couples the anti-self-duality 
equation on $\R\times[0,1]\times\S$ with Lagrangian boundary conditions.
We expect that two different degenerations of the metric on $[0,1]\times\S$ 
will give rise to isomorphisms
that would prove the Atiyah-Floer conjecture~\footnote{
There are moreover product structures on all three Floer homologies that 
should be intertwined by the isomorphisms, as sketched in \cite{Sa}. 
Our analytic setup should allow for their definition and identification, 
but we do not discuss this topic here.
}
\begin{align} \label{isom1}
\HF_*^{\rm inst}([0,1]\times\S,\cL_{H_0}\times\cL_{H_1})
&\cong 
\HF_*^{\rm inst}(H_0\cup_\S H_1) , \\
\label{isom2}
\HF_*^{\rm inst}([0,1]\times\S,\cL_{H_0}\times\cL_{H_1})
&\cong 
\HF_*^{\rm symp}(\cR_\S,L_{H_0},L_{H_1}) .
\end{align}
This approach separates the difficulties:
The first isomorphism is a purely gauge theoretic comparison between anti-self-dual
instantons over domains with and without boundary.
The second isomorphism requires a comparison between anti-self-dual instantons and 
pseudoholomorphic curves (both with Lagrangian boundary conditions), that would 
be a generalization of the adiabatic limit of Dostoglou-Salamon \cite{DS}, 
which they used to prove an analogon of the Atiyah-Floer conjecture for mapping tori.
The mapping torus case does not involve boundary conditions. Moreover, the 
underlying bundle is nontrivial so that the moduli space of flat connections is smooth.
In contrast, the Heegard splitting case deals with trivial bundles for which
the moduli space $\cR_\S$ and its Lagrangian submanifolds are always singular.

So the Atiyah-Floer conjecture poses as a first task (which we do not approach here)
the construction of a symplectic Floer homology for symplectic and Lagrangian manifolds 
with quotient singularities.
In fact, the singular symplectic space $\cR_\S$ is the symplectic quotient 
(in the sense of Atiyah and Bott \cite{AB}) of a Hamiltonian group action 
(the infinite dimensional gauge group) on an infinite dimensional symplectic space 
(the space of connections over a Riemann surface).
In the case of a finite dimensional Hamiltonian group action with 
smooth and monotone symplectic quotient, Gaio and Salamon \cite{GS} have identified 
the Gromov-Witten invariants of the symplectic quotient with new invariants arising 
from the symplectic vortex equations.

The anti-self-duality equation on $\R\times[0,1]\times\S$ 
is the exact analogue of the symplectic vortex equations for $\cR_\S$.
We will show in section~\ref{sec:FH} that the analytic behaviour of these trajectories
of the new Floer homology is a mixture of local effects in the interior 
-- as they are expected for anti-self-dual instantons --
and surprising semiglobal effects near the boundary that resemble to the 
behaviour of pseudoholomorphic curves in $\cR_\S$.
This shows that the new Floer homology indeed provides a good interpolation between
the two Floer homologies in the Atiyah-Floer conjecture.\\

More generally, an instanton type Floer homology for $3$-manifolds with boundary 
should naturally use Lagrangian boundary conditions.
Fukaya \cite{Fu} gives such a setup in the case of a nontrivial bundle:
The anti-self-duality equation is coupled via a degeneration of the metric
to the pseudoholomorphic curve equation in the moduli space of flat connections
(which is smooth in this case).
Our new trajectory equation is a different setting that arises naturally from the 
Chern-Simons functional -- the Morse function in the instanton Floer theory.
It works 
in the gauge theoretic setting up 
to the boundary, which has the advantage that the Lagrangians are smooth Banach 
submanifolds of a symplectic Banach space, although the quotients might be singular.
We thus give a setup for an instanton Floer homology $\HF_*^{\rm inst}(Y,\cL)$
of a compact 3-manifold $Y$ with boundary and a gauge invariant Lagrangian submanifold
$\cL$ in the space of $\SU(2)$-connections over $\pd Y$.

This program is carried through in \cite{SW} for the case where 
$\cL=\cL_H$ arises from a disjoint union of handle bodies $H$ with boundary 
$\pd H=\pd Y=\S$ such that $Y\cup_\S H$ is a homology $3$-sphere.
We expect that the isomorphism (\ref{isom1}) will be true in this more general
setting, 
\begin{equation} \label{isom1'} \tag{1'}
\HF_*^{\rm inst}(Y,\cL_H) \cong \HF_*^{\rm inst}(Y\cup_\S H) .
\end{equation}
The assumption $\cL=\cL_H$ is more of technical nature and is not required for
the basic compactness in theorem~\ref{thm:C}.
We also have an approach to removing this assumption in theorems~\ref{thm:EQ} 
and \ref{thm:RemSing}, 
based on Mrowka's understanding of the gauge group in borderline Sobolev cases. 
More essentially, we need the nontrivial flat connections on $Y$ with boundary 
condition in $\cL$ to be irreducible (i.e.\ to have discrete isotropy in the gauge group).
For the same reason, the original instanton Floer homology 
is only defined for homology $3$-spheres.
Now starting from a Heegard splitting $H_0\cup_\S H_1$ of a general closed $3$-manifold,
the irreducibility could be achieved by perturbing the Lagrangians $\cL_{H_0}$ and $\cL_{H_1}$.
Thus the problem of reducible connections can be transferred to transversality 
questions in our new instanton Floer homology with Lagrangian boundary conditions.\\

Section~\ref{sec:gauge} provides an introduction to the gauge theoretic background.
For the symplectic background we refer to~\cite{MS}.
We explain the Chern-Simons functional and 
the moment map picture of the gauge group action
and give the setup for an instanton Floer homology 
$\HF_*^{\rm inst}(Y,\cL)$ with Lagrangian boundary conditions.
In section~\ref{sec:FH} we specialize to the case $Y=[0,1]\times\S$ and 
the Lagrangian submanifold $\cL_H=\cL_{H_0}\times\cL_{H_1}$ arising from two 
handle bodies ${H=H_0\sqcup H_1}$ such that $H_0\cup_\S H_1\cong Y\cup_{\S\sqcup\S} H$ 
is a homology $3$-sphere.
We give a detailed account of the new Floer homology 
${\HF_*^{\rm inst}([0,1]\times\S,\cL_{H_0}\times\cL_{H_1})}$,
comparing its definition and the analytic properties of its trajectories
to those of $\HF_*^{\rm inst}(H_0\cup_\S H_1)$ 
and $\HF_*^{\rm symp}(\cR_\S,L_{H_0},L_{H_1})$ (what it would be if 
these quotients were smooth).
In section~\ref{sec:AF} we sketch the ideas for proofs of the isomorphisms 
(\ref{isom1}) and~(\ref{isom2}).

The last two sections are a rough guide to the analysis of anti-self-dual 
instantons with Lagrangian boundary conditions, which was established in 
\cite{W Cauchy, W elliptic, W bubbling, W mean} in full technical detail.
Section~\ref{sec:L} provides an overview of the properties of gauge invariant 
Lagrangian submanifolds in the space of connections over a Riemann surface.
It moreover describes the special extension properties of Lagrangian submanifolds
that arise from handle bodies. 
In section~\ref{sec:analysis} we sketch the proofs of the analytic results 
in section~\ref{sec:FH}.
We put the proofs into context with the standard proofs
of Uhlenbeck compactness (for anti-self-dual instantons)
and Gromov compactness (for pseudoholomorphic curves)
since -- just as the results -- each proof requires a subtle combination
of the best techniques from both gauge theory and symplectic topology, which 
we hope the reader will find entertaining.

\pagebreak

\section{Gauge theory and symplectic topology}
\label{sec:gauge}

We give an introduction to some gauge theoretic concepts and notations.
More details and proofs can be found in e.g.\ \cite{DK,W}.

Let $\rG$ be a compact Lie group. The Lie algebra $\cg=\rT_{\smone}\rG$ 
is equipped with a Lie bracket $[\cdot,\cdot]$ and with a
$\rG$-invariant inner product $\la\cdot,\cdot\ra$.
For the instanton Floer theories we will be using $\rG=\SU(2)$
with the commutator $[\x,\e]=\x\e - \e\x$ and the trace
$\la\x,\e\ra =-\tr(\x\e)$ for ${\x,\e\in\su(2)}$.
We describe a {\bf connection} on the trivial $\rG$-bundle 
$\rG\times X \to X$ over a manifold $X$
as a $\cg$-valued $1$-form ${A\in\Om^1(X;\cg)}$ 
and thus denote the space of smooth connections by 
$$
\cA(X):=\Om^1(X;\cg).
$$
(The discussion in this section generalizes to nontrivial bundles, 
where connections are given by $1$-forms with values in an associated bundle.)
On the trivial bundle a $1$-form $A\in\cA(X)$ 
corresponds to an equivariant distribution
${\{(-gA(Y),Y)\st Y\in\rT_xX\}\subset\rT_{(g,x)} (\rG\times X)}$
 of horizontal subspaces.
The corresponding covariant derivative 
on sections $s:X\to E$ of a trivial vector bundle 
with structure group $\rG\subset\Hom(E)$
is $\nabla_A s : Y \mapsto \nabla s(Y) + A(Y) s$.


The {\bf curvature} of a connection $A\in\cA(X)$ is given by the $2$-form
$$
F_A :=\rd A + \half [A\wedge A]\;\in\Om^2(X;\cg) .
$$
Throughout $[\cdot\wedge\cdot]$ indicates that the values of the differential forms 
are paired by the Lie bracket. 
The differential of the map $A\mapsto F_A$ at a connection $A\in\cA(X)$
is the 'twisted' exterior derivative $\rd_A:\Om^1(X;\cg)\to\Om^2(X;\cg)$. 
In general, $\rd_A:\Om^k(X;\cg)\to\Om^{k+1}(X;\cg)$ 
acts on $\cg$-valued differential forms by
$$
\rd_A\e := \rd\e + [A\wedge\e] .
$$
One checks that $\rd_A\rd_A\e = [F_A\wedge\e]$, so ${\rd_A}^2=0$ 
iff the curvature vanishes.
Such connections are called {\bf flat} and we denote the set of flat connections~by
$$
\Af(X):=\{ A\in\cA(X) \st F_A=0 \} .
$$
Moreover, a connection is flat iff the
horizontal distribution is locally integrable.
So parallel transport with respect to a flat connection
around a loop is given by an element in the group $\rG$ that is 
invariant under homotopy of the loop with fixed base point $x\in X$.
Thus the {\bf holonomy} induces a map
$$
\hol_x : \Af(X) \to \Hom(\pi_1(X,x),\rG) .
$$
Next, connections that are the same up to a bundle isomorphism are called gauge equivalent.
The bundle isomorphisms of the trivial bundle can be identified with maps $u:X\to\rG$ that are
called {\bf gauge transformations}.
Composition of bundle isomorphisms corresponds to multiplication of gauge transformations,
so the space of smooth gauge transformations has the structure of a group, called the {\bf gauge group}
$$
\cG(X):=\cC^\infty(X,\rG) .
$$
The action of $\cG(X)$ on the space of connections $\cA(X)$, called the {\bf gauge action}, 
is given by pullback of the connection 
(i.e.\ the horizontal subspace or the covariant derivative), hence
$$
u^*A := u^{-1}A u + u^{-1}\rd u \qquad\text{for}\; u\in\cG(X), A\in\cA(X) .
$$
The space of flat connections $\Af(X)$ is obviously invariant under $\cG(X)$, 
and the curvature transforms by $F_{u^*A}=u^{-1} F_A u$.
The holonomy of a flat connection $A\in\Af(X)$ transforms by conjugation under
$u\in\cG(X)$, more precisely $\hol_x(u^*A)=u(x)^{-1}\hol_x(A) u(x)$ for the holonomy
based at ${x\in X}$. Similarly, a change of the base point also transforms the 
holonomy by conjugation. Hence the holonomy descends to a map
$$
\hol : \Af(X)/\cG(X) \to \Hom(\pi_1(X),\rG)/\rG =: \cR_X  ,
$$
where the action of $\rG$ is by conjugation.
If there are no nontrivial $\rG$-bundles over $X$,\footnote
{
This is for example the case for $\dim X =2$ or $3$ and a connected, 
simply connected group as $\rG=\SU(2)$.
It also holds for a handle body $X$ and any connected group $\rG$.
}
then this is in fact an isomorphism and we will identify the representation space 
$\cR_X$ with $\Af(X)/\cG(X)$.
In general this is an isomorphism when 
taking the union over all isomorphism classes of bundles 
on the left hand side, see e.g.\ \cite[Proposition 2.2.3]{DK}.

\subsection*{Uhlenbeck compactness} \hspace{2mm} 

The observations above shows that the 
moduli space of flat connections is a compact subset of $\cA(X)/\cG(X)$ 
(in the $\cC^\infty$-topology).
Uhlenbeck's weak compactness theorem is a remarkable generalization 
of this compactness to connections with small curvature.
It is the starting point for all analysis in gauge theory, so this is a good 
point to introduce the Sobolev completions of the spaces of connections and 
gauge transformations. For a compact manifold $X$ and 
for $k\in\N_0$ and $1\leq p\leq\infty$ let
\begin{align*}
\cA^{k,p}(X) := W^{k,p}(X,\rT^*X\otimes\cg), \qquad
\cG^{k,p}(X) := W^{k,p}(X,\rG).
\end{align*}
For $kp>\dim X$ the gauge group $\cG^{k,p}(X)$ is a Banach manifold, on which 
multiplication and inversion are smooth, and it acts smoothly on $\cA^{k-1,p}(X)$.
We equip $X$ with a metric, then for any $p\geq 1$ the $L^p$-norm of the curvature,
$$
\| F_A \|_p^p = \int_X |F_A|^p ,
$$
is a gauge invariant quantity.
For $p=\half\dim X$ this is the conformally invariant Yang-Mills energy of the connection,
which can concentrate at single points.
Thus for ${p\leq \half\dim X}$ one cannot expect the compactness of a set of connections with 
bounded $L^p$-norm of the curvature.
Uhlenbeck's result \cite{U2} says that for $p>\half\dim X$ 
however, every such set is compact in the weak
$W^{1,p}$-topology on the quotient $\cA^{1,p}(X)/\cG^{2,p}(X)$.

\begin{thm} {\bf (Weak Uhlenbeck Compactness)} \label{thm:weak comp}\\
Let $X$ be a compact manifold and let $p>\half\dim X$.
Suppose that $A_i\in \cA(X)$ is a sequence of connections
such that $\|F_{A_i}\|_p$ is uniformly bounded.
Then, after going to a subsequence, there exists a 
sequence of gauge transformations $u_i\in\cG(X)$ such that 
$u_i^*A_i \to A_\infty$ converges in the weak $W^{1,p}$-topology
to a connection $A_\infty\in\cA^{1,p}(X)$.
\end{thm}

In fact, one even has a weak $W^{1, \dim X/2}$-compactness 
if one assumes that 
every point in $X$ has a neighbourhood on which the 
Yang-Mills energy is bounded by a small constant.
We will need the slightly stronger $W^{1,p}$-compactness since it allows us 
to globally (not just locally over small balls in $X$)
work in in a local slice of the gauge action.
The local slice theorem says that any connection $A'$ that is suitably close 
to a fixed reference connection $A$ can be put into relative Coulomb gauge,
i.e.\ $u^*A'-A$ is $L^2$-orthogonal to the gauge orbit through~$A$.
The linearized gauge action $\rT_\smone\cG(X)\to\rT_A\cA(X)$ at $A\in\cA(X)$
is given by $\rd_A:\Om^0(X;\cg)\to\Om^1(X;\cg)$.
Its formal adjoint is $\rd_A^*=-*\rd_A*$.
More generally, with $m(X,k)=(\dim X - k)(k-1)$ the twisted coderivative is
$$
\rd_A^* := - (-1)^{m(X,k)} * \rd_A * \;:\;\Om^k(X;\cg)\to\Om^{k-1}(X;\cg).
$$
Here we only give the sequential form of the local slice theorem. 
A stronger statement and proof can be found e.g.\ in \cite[Theorem~F]{W}.

\begin{thm} {\bf (Local Slice Theorem)} \label{thm:local slice}\\
Let $X$ be a compact Riemannian manifold with smooth boundary and
let $p>\half\dim X$.
Suppose that $A_i\in\cA^{1,p}(X)$ is a sequence of connections
such that $A_i\to A\in\cA(X)$ in the weak $W^{1,p}$-topology.
Then for sufficiently large $i$ 
there exist gauge transformations $u_i\in\cG^{2,p}(X)$ such that
$u_i^*A_i\to A$ and
\[
\left\{\begin{aligned}
\rd_A^*(u_i^*A_i-A)&=0 , \\
*(u_i^*A_i-A)|_{\pd X}&=0 .
\end{aligned}\right.
\]
\end{thm}

To see the strength of these two theorems consider the following example:
Let $X$ be a compact Riemannian $4$-manifold.
The extrema of the Yang-Mills energy $\int_X |F_A|^2$ are called
{\bf Yang-Mills instantons}. They satisfy the equation $\rd_A^*F_A=0$
and the boundary condition $*F_A|_{\pd X}=0$, which 
-- augmented with the local slice conditions --
pose an elliptic boundary value problem. 
Uhlenbeck's compactness theorem for Yang-Mills instantons with $L^p$-bounded curvature
then is a corollary of theorems~\ref{thm:weak comp} and \ref{thm:local slice}.

\begin{thm}\rm 
{(\bf Strong Uhlenbeck Compactness): }\\
Let $A_i\in \cA(X)$ be a sequence of Yang-Mills instantons
such that $\|F_{A_i}\|_p$ is uniformly bounded for some $p>2$.
Then, after going to a subsequence, there exists a 
sequence of gauge transformations $u_i\in\cG(X)$ such that 
$u_i^*A_i \to A_\infty$ converges in the $\cC^\infty$-topology
to another Yang-Mills instanton $A_\infty\in\cA(X)$.
\end{thm}

{\bf Anti-self-dual instantons} on an oriented Riemannian $4$-manifold $X$ 
are solutions $A\in\cA(X)$ of the first order equation
$$
F_A+*F_A=0.
$$
By the Bianchi identity $\rd_A F_A=0$ these are special
solutions of the Yang-Mills equation $\rd_A^*F_A=0$.
On a manifold with boundary however, 
the anti-self-duality equation with boundary condition $*F_A|_{\pd X}=0$ 
is an overdetermined boundary value problem comparable to Dirichlet boundary 
conditions for holomorphic maps.
This is another reason why it is natural to consider (weaker) Lagrangian boundary conditions
for anti-self-dual instantons.

\subsection*{The moduli space of flat connections over a Riemann surface} \hspace{2mm} 

Let $\S$ be a Riemann surface.
The natural symplectic form on the space of connections $\cA(\S)=\Om^1(\S;\cg)$ is
\begin{equation} \label{omega}
\o(\a,\b) := \int_\Si \la \a \wedge \b \ra 
\qquad\quad\text{for}\;\; \a,\b\in\Om^1(\S;\cg) .
\end{equation}
Here and throughout $\la \cdot \wedge \cdot \ra$ indicates that the values 
of the differential forms are paired by the inner product. 
Note that for any metric on $\Sigma$ the Hodge operator $*$ is a complex structure
on $\cA(\S)$, which is compatible with $\o$ and induces the $L^2$-metric
$\o(\a,*\b)=\la \a,\b\ra_{L^2}$.

It was observed by Atiyah and Bott \cite{AB} that the action of the gauge group 
$\cG(\S)$ on $\cA(\S)$ can be viewed as Hamiltonian action of an infinite dimensional Lie group.
The Lie algebra of $\cG(\S)$ is $\Om^0(\Si;\cg)$ and the infinitesimal action of
$\x\in\Om^0(\Si;\cg)$ is given by the vector field
$$
X_\x : \cA(\S) \rightarrow  \Om^1(\S;\cg) , \quad
X_\x(A) = \rd_A\x . 
$$
This is the Hamiltonian vector field of the function 
$A \mapsto \int_\S\la \m(A) \,,\, \x \ra$, where
$$
\m :\cA(\S) \rightarrow \Om^0(\S;\cg) , \quad
\m(A) = *F_A 
$$
can be considered as a moment map.
Its differential is $\rd\mu(A)=*\rd_A$, 
so one indeed has for all $\b\in \Om^1(\S;\cg)$
$$
\o( X_\x(A) \,,\,\b ) \;=\; \int_\Si \la \rd_A\x \wedge \b \ra 
\;=\; - \int_\S \la \x \,,\, *\rd_A\b \ra 
\;=\; - \int_\S\la \rd\mu(A) \b \,,\, \x \ra.
$$
The zero set of $\m$ is the set of flat connections.
So the moduli space of flat connections on $\S$ can be seen as the
symplectic quotient of the gauge action, 
$$
\cR_\S \;=\; \cA_{\rm flat}(\S)/\cG(\S)
\;=\; \m^{-1}(0)/\cG(\S) 
\;=\; \cA(\S) /\hspace{-1mm}/ \cG(\S) .
$$
This quotient $\cR_\S \cong \Hom(\pi_1(\S),\rG)/\rG$
is singular at the reducible representations, but 
for irreducible\footnote{
A connection $A\in\cA_{\rm flat}(\S)$ is called irreducible if its isotropy subgroup
of $\cG(\S)$ (the group of gauge transformations that leave $A$ fixed) is discrete,
i.e.\ $\rd_A|_{\Om^0}$ is injective. For a closed Riemann surface this is equivalent to
$\rd_A^*|_{\Om^1}$ being surjective.}
$A\in\Af(\S)$ 
it is a smooth manifold near the gauge equivalence class $[A]$.
To understand its tangent space, notice that
the linearized action $\rT_\smone\cG(\S) \to \rT_A\cA(\S)$ is
given by $\rd_A:\Om^0(\Si;\cg) \to \Om^1(\Si;\cg)$.
At a flat connection $A\in\Af(\S)$ it fits into a chain complex with 
the differential of the moment map (since $\rd_A\comp\rd_A=0$),
$$
\Om^0(\Si;\cg) 
\overset{\rd_A}{\longrightarrow}  \Om^1(\Si;\cg)  
\overset{*\rd_A}{\longrightarrow}  \Om^0(\Si;\cg) .
$$
So the tangent space to $\cR_\S$ at $[A]$ is the twisted first homology group, which can 
be identified with the harmonic $1$-forms,
$$
\rT_{[A]}\cR_\S \;=\; \ker*\rd_A / \im\rd_A 
\;\cong\; \ker\rd_A\cap\ker\rd_A^* \;=:\, h^1_A .
$$
Hodge theory gives a corresponding $L^2$-orthogonal splitting
\begin{equation}\label{Hodge}
\Om^1(\S;\cg)=\im\rd_A \oplus \im(*\rd_A) \oplus h^1_A ,
\end{equation}
where $\im\rd_A\oplus h^1_A=\ker*\rd_A$ is the tangent space to $\m^{-1}(0)=\Af(\S)$
and $\im(*\rd_A)\oplus h^1_A=(\im\rd_A)^\perp$ is the local slice
of the gauge action through $A$.

We have seen that the moduli space of flat connections $\cR_\S$ is a smooth manifold 
of dimension $(2g-2)\dim\rG$ with singularities at the reducible connections. 
Moreover, the symplectic structure (\ref{omega}) on $\cA(\S)$ is $\cG(\S)$-invariant
and induces a symplectic structure 
on the smooth part of $\cR_\S$. For harmonic representatives 
$\a,\b\in h^1_A\cong\rT_{[A]}\cR_\S$ it is again given by
$\o(\a,\b)=\int_\S \la \a\wedge\b \ra$.
In this representation of the tangent space we also see that the 
Hodge operator $*$ descends to $\cR_\S$. 
So $(\cR_\S,\o)$ is a (singular) symplectic manifold with compatible 
almost complex structure~$*$.

\subsection*{The Chern-Simons functional and instanton Floer homology}  \hspace{2mm} 

Let $Y$ be a compact oriented $3$-manifold. The {\bf Chern-Simons 1-form}
$\l$ on the space of connections $\cA(Y)$ is given by
$$
\l_A (\a) :=  \int_Y \la F_A \wedge \a \ra 
\qquad\quad\text{for}\;\; \alpha\in\rT_A\cA(Y)=\Om^1(Y;\cg).
$$
This $1$-form is equivariant, $\l_{u^*A}(u^{-1}\a u)=\l_A(\a)$.
If $\pd Y=\emptyset$ or ${F_A|_{\pd Y}=0}$, then $\l$ is also horizontal and thus 
descends to the (singular) moduli space\footnote{
Note that $\cB(Y)$ is not the moduli space of flat connections $\cR_Y$, 
but the infinite dimensional and singular space of all connections 
modulo gauge equivalence.
}${\cB(Y):=\cA(Y)/\cG(Y)}$.
Indeed, a tangent vector to the gauge orbit through ${A\in\cA(Y)}$ has the form
$\a=\rd_A\x$ with $\x\in\Om^0(Y;\cg)$, and by Stokes' theorem and the Bianchi identity
\begin{equation}\label{hor}
\l_A(\rd_A\x) \;=\;
- \int_Y \la \rd_A F_A \,,\, \x \ra 
+ \int_{\pd Y} \la F_A \,,\, \x \ra 
\;=\; 0 .
\end{equation}
To calculate the differential of $\l$ consider $\a,\b\in\Om^1(Y;\cg)$ as (constant) 
vector fields on $\cA(Y)$, then their Lie bracket vanishes and
\begin{align}
\rd\l(\a,\b) 
&\;=\; \nabla_\a (\l(\b)) - \nabla_\b (\l(\a)) \nonumber\\
&\;=\; \int_Y \la \rd_A\a\wedge\b \ra - \int_Y \la \rd_A\b\wedge\a \ra
\;=\; \int_{\pd Y} \la \a\wedge\b \ra. \label{closed}
\end{align}
So for $\pd Y=\emptyset$ the Chern-Simons 1-form descends to a closed
$1$-form on $\cB(Y)$.
In fact, $\l$ is the differential of the {\bf Chern-Simons functional}
$$
\CS(A) := \half \int_Y \la A \wedge \bigl( F_A - \tfrac 16 [A\wedge A] \bigr) \ra  .
$$
For a more illuminating definition let $X$ be a compact $4$-manifold with 
boundary $\pd X=Y$, then for any $\tA\in\cA(X)$ with $\tA|_{\pd X}=A$
$$
\CS(A) = \half \int_X \la F_\tA \wedge F_\tA \ra  .
$$
For closed $X$ the right hand side is a topological invariant of the bundle. 
We fix $\rG=\SU(2)$, then this invariant is $4\pi^2 c_2(P)$,
where $\tA$ is a connection on the bundle $P\to X$.
From this one can see that the Chern-Simons functional 
descends to an $S^1$-valued functional $\CS:\cB(Y)\to\R/4\pi^2\Z$ 
since it changes by $\CS(A) - \CS(u^*A) = 4\pi^2\deg(u)\in 4\pi^2\Z$
under gauge transformations.
If $Y$ is a homology $3$-sphere, then Floer \cite{F1} used the generalized Morse 
theory for this functional to define the instanton Floer homology $\HF^{\rm inst}_*(Y)$.
Roughly speaking, the Floer complex is generated by the zeros of $\rd\CS=\l$, 
i.e.\ by the flat connections $A\in\Af(Y)$ modulo $\cG(Y)$.
The differential on the complex is defined by counting negative gradient flow lines, 
so we choose a metric on $Y$ and thus fix an $L^2$-metric on $\cA(Y)$.
Then the gradient of $\CS$ is $A\mapsto *F_A$
and a negative gradient flow line is a path $A:\R\to\cA(Y)$ satisfying
$$
\pd_s A = - * F_A  .
$$
Equivalently, one can view this path as connection $\X= \P\ds + A \in\cA(\R\times Y)$
in the special gauge $\P\equiv 0$. 
Then the above equation is the anti-self-duality
equation $F_\X+*F_\X=0$ for $\X$. 
For a general connection $\X\in\cA(\R\times Y)$ this so-called 
temporal gauge can always be
achieved by the solution $u\in\cG(\R\times Y)$ of ${\pd_s u = -\P u}$ with
$u|_{s=0}\equiv\one$.
So the negative gradient flow lines of the Chern-Simons functional modulo $\cG(Y)$
are in one-to-one correspondence with the anti-self-dual connections
$\X\in\cA(\R\times Y)$ modulo $\cG(\R\times Y)$. 
An extensive discussion of instanton Floer homology for closed $3$-manifolds can be
found in Donaldson's book \cite{Donaldson book}.\\

If $Y$ has nonempty boundary $\pd Y=\Si$, then the differential (\ref{closed}) 
is the symplectic form $\o$ on $\a|_\S,\b|_\S\in\cA(\Si)$, compare (\ref{omega}).
To render $\l$ closed, it is natural\footnote{
If $\cL\subset\cA(\Si)$ is any submanifold, then the closedness of $\l|_{\cA(Y,\cL)}$
is equivalent to $\o|_\cL\equiv 0$, and the maximal such submanifolds are precisely
the Lagrangian submanifolds.
}
to pick a Lagrangian submanifold 
$\cL\subset\cA(\S)$ and restrict $\l$ to 
$$
\cA(Y,\cL):=\{A\in\cA(Y)\,|\,A|_\S\in\cL\} .
$$
More precisely, we fix a $p>2$ 
and make the following assumptions to ensure that 
$\l$ defines a closed 1-form on $\cB(Y,\cL):=\cA(Y,\cL)/\cG(Y)$.
\begin{enumerate} 
\item 
\label{Lag ass}
$\cL\subset\cA^{0,p}(\S)$ is a Banach submanifold that is isotropic, $\o|_\cL\equiv 0$,
and coisotropic in the sense of the following implication for all $\a\in\cA^{0,p}(\S)$:
If $\o(\a,\b)=0$ for all $\b\in\rT_A\cL$, then $\a\in\rT_A\cL$.
\item $\cL$ is invariant under $\cG^{1,p}(\S)$.
\item
$\cL\subset\cA_{\rm flat}^{0,p}(\S)$ lies in the space of weakly flat 
connections.\footnote{
See~\cite[Sec.~3]{W Cauchy} or use the fact
$\cA_{\rm flat}^{0,p}(\S)=\cG^{1,p}(\S)^*\Af(\S)$ as definition.
The conditions (ii) and (iii) are equivalent when $\rG$ is connected and 
$[\cg,\cg]=\cg$, e.g.\ for $\SU(2)$.
}
\end{enumerate}
Here (ii) ensures that $\cG(Y)$ acts on $\cA(Y,\cL)$, and (iii) implies that
$\l$ is horizontal by (\ref{hor}).
These assumptions also imply that $\cL$ descends to a 
(singular) Lagrangian submanifold in the 
(singular) moduli space of flat connections,
$$
L := \cL/\cG^{1,p}(\Si) \subset \cR_\Si = \cA_{\rm flat}(\Si)/\cG(\Si) .
$$
The assumptions (i)-(iii) also imply the orthogonal splitting, see
section~\ref{sec:L},
$$
\Om^1(\S;\cg)=\rT_A\cL \;\oplus *\rT_A\cL \qquad\text{for all}\; A\in\cL .
$$
Compare this to (\ref{Hodge}) and note that $\im\rd_A\subset\rT_A\cL \subset \ker\rd_A$ 
due to (ii),(iii). So the $\rT_A\cL$ are determined up to a choice of Lagrangian
subspaces in $h^1_A$.
Conversely, any Lagrangian $L\subset\cR_\S$ lifts to a (possibly nonsmooth)
${\cL\subset\cA^{0,p}(\S)}$ as above.
In order to obtain a well defined Floer homology one should moreover
assume that $L$ is simply connected (which ensures a monotonicity property).
In general, $\cL$ is not simply connected, but its fundamental group cancels with that 
of $\cG(\Si)$. This is the reason why $\l$ is not exact but can only be written
as the differential of the multi-valued Chern-Simons functional
$$
\CS_\cL(A) 
= \half \int_Y \la A\wedge \bigl( F_A - \tfrac 16 [A\wedge A]\bigr) \ra 
 +\int_0^1\int_\Si \la \tA(t)\wedge \pd_t\tA(t)\ra \,\dt .
$$
This involves the choice of a path $\tA:[0,1]\to\cL$ with
$\tA(1)=A|_\Si$ and ${\tA(0)=A_0}$ a fixed reference connection in $\cL$.
Again we fix $\rG=\SU(2)$ to obtain a functional
$\CS_\cL : \cA(Y,\cL) \to \R/{4\pi^2\Z}$, 
which directly descends to $\cB(Y,\cL)$ since
the gauge group $\cG(Y)$ is connected.
We now propose to define a new Floer homology $\HF^{\rm inst}_*(Y,\cL)$ from the 
generalized Morse theory of the functional $\CS_\cL : \cB(Y,\cL) \to \R/{4\pi^2\Z}$.

A critical point in this theory is a flat connection $A\in\Af(Y)$ with Lagrangian
boundary condition $A|_\Si\in\cL$ (modulo $\cG(Y)$), and a negative gradient flow line 
is a path $A:\R\to\cA(Y)$ (modulo $\cG(Y)$) satisfying   
\begin{equation} \label{new crit}
\pd_s A = - * F_A , \qquad A(s)|_\Si\in\cL \quad\forall s\in\R .
\end{equation}
Again, this is the anti-self-duality equation for 
$\X= A + \P\ds \in\cA(\R\times Y)$ in the temporal gauge $\P\equiv 0$. 
So the gauge equivalence classes of the gradient flow lines 
are in one-to-one correspondence with the gauge equivalence classes of
anti-self-dual instantons with Lagrangian boundary conditions, i.e.\
solutions $\X\in\cA(\R\times Y)$ of the boundary value problem
\begin{equation} \label{new flow}
 F_\X + * F_\X = 0   ,\qquad
 \X|_{\{s\}\times\Si} \in\cL\qquad \forall s\in\R.
\end{equation}

\subsection*{Lagrangians and handle bodies} \hspace{2mm}

We have seen before how a Riemann surface $\S$ gives rise to a 
(singular) symplectic manifold $\cR_\S=\Hom(\pi_1(\S),\rG)/\rG$,
which is a finite dimensional reduction of a symplectic Banach space
$\cA(\S)=\Om^1(\S;\cg)$ that arises from gauge theory on $\S$.
We will now discuss a class of examples of Lagrangian Banach submanifolds
$\cL_H\subset\cA(\S)$ that arise from gauge theory on a handle body $H$ with 
$\pd H=\S$, and that reduce to finite dimensional (singular) Lagrangian
submanifolds $L_H\cong\Hom(\pi_1(H),\rG)/\rG\subset\cR_\S$.
Here and throughout a handle body is an oriented $3$-manifold with boundary 
that is obtained from the $3$-ball by attaching a finite number of $1$-handles.

For this purpose let $\rG$ be a compact, connected,
and simply connected Lie group (e.g.\ $\rG=\SU(2)$) and
let $\S$ be a Riemann surface.
For a start let $H$ be any compact $3$-manifold with boundary $\pd H=\S$. Then
$$
\cL_H:=\bigl\{ \tA|_\S \st \tA\in\Af(H) \bigr\}
\subset \cA(\S)
$$
satisfies the assumptions (ii) $\cL_H\subset\Af(\S)$,
(iii) $\cG(\S)^*\cL_H=\cL_H$,\footnote{
This uses the assumption 
$\pi_1(\rG)=\{0\}$ and the fact that $\pi_2(\rG)=\{0\}$ for every
compact Lie group, so that every gauge transformation $u:\S\to\rG$ 
extends to $\tu:H\to\rG$.}
and is isotropic: Consider paths
$\tA_1,\tA_2:(-\ep,\ep)\to\Af(H)$ with $\tA_i(0)=\tA$,
then ${\rd_\tA\pd_t\tA_i(0)=\pd_t\bigr|_{t=0}F_{\tA_i}=0}$ and hence
with the symplectic form (\ref{omega})
\begin{align*}
&\o\bigl(\pd_t\tA_1(0) \,,\, \pd_t\tA_2(0)\bigr) 
\;=\;\tint_{\pd H} \la \pd_t\tA_1(0) \wedge \pd_t\tA_2(0) \ra \\
&=\tint_H \la \rd_\tA\pd_t\tA_1(0) \wedge \pd_t\tA_2(0) \ra
           - \la \pd_t\tA_1(0) \wedge \rd_\tA\pd_t\tA_2(0) \ra 
\;=\;0.
\end{align*}
So $\cL_H$ descends to an isotropic subset in the symplectic quotient
$$
L_H \,:=\; \cL_H/\cG(\S) \;\subset\; \cR_\S \;=\; \cA(\S)/\hspace{-1mm}/\cG(\S) .
$$
The holonomy provides an isomorphism
$$
L_H \,\cong\; \Hom\bigl(\tfrac{\pi_1(\S)}{\pd\pi_2(H,\S)},\rG\bigr)/\rG 
\;\subset\; \cR_\S \;\cong\; \Hom(\pi_1(\S),\rG)/\rG .
$$
This is since the holonomy of a flat connection on $H$ is trivial
on the contractible loops in $\pd\pi_2(H,\S)$; and all representations
of $\pi_1(\S)/\pd\pi_2(H,\S)$ can be realized by a flat connection on $H$ 
since that quotient embeds into $\pi_1(H)$ by the long exact sequence for 
homotopy
$$
\dots \to \pi_2(H)
\to \pi_2(H,\S)
\overset{\pd}{\to} \pi_1(\S)
\overset{\i}{\to} \pi_1(H)
\to \pi_1(H,\S)
\to \dots
$$
This also shows that 
$\frac{\pi_1(\S)}{\pd\pi_2(H,\S)}\cong\pi_1(H)$ 
if $H$ is a handle body (so $\pi_2(H)$ and $\pi_1(H,\S)$ vanish).
Now consider the commuting diagram of long exact sequences
for homology and cohomology with the vertical Poincare duality:
$$
\begin{array}{ccccc}
H_2(H,\S) &
\overset{\pd}{\longrightarrow} & H_1(\S) &
\overset{\i}{\longrightarrow} & H_1(H)  \\
\updownarrow & &\updownarrow& &\updownarrow\\
H^1(H) &
\overset{\i^*}{\longrightarrow} &  H^1(\S)  &
\overset{\pd^*}{\longrightarrow} &  H^2(H,\S)
\end{array}
$$
One can read off that
${(\im\pd)^\perp \cong (\im\i^*)^\perp = (\ker\pd^*)^\perp = \im\pd}$,
and we obtain ${\dim\frac{H_1(\S)}{\pd H_2(H,\S)}=\half\dim H_1(\S)}$.
Hence $\dim L_H = \half \dim\cR_\S$ at smooth points.
So a general compact $3$-manifold $H$ with $\pd H=\S$ gives rise to a 
(singular) Lagrangian $\cL_H\subset\cR_\S$,
and in fact $\cL_H\subset\cA(\S)$ is Lagrangian up to possible singularities.
If $H$ is a handle body, then one can prove that $\cL_H$ is in fact smooth, \footnote{
Its $L^p$-completion is a Banach submanifold of $\cA^{0,p}(\S)$, see section~\ref{sec:L}.}
which is essentially due to the fact that
$\frac{\pi_1(\S)}{\pd\pi_2(H,\S)}\cong\pi_1(H)$ is a free group.

This correspondence between low dimensional topology, symplectic topology, 
and gauge theory is summarized in a table on page~\pageref{table 1}.
To round off this discussion, note that a Heegard splitting $H_0\cup_\S H_1$
of a $3$-manifold into two handle bodies $H_0, H_1$ with common boundary
$\pd H_i=\S$ gives rise to a pair of (singular) Lagrangians in a 
symplectic manifold, $L_{H_0},L_{H_1}\subset\cR_\S$.
Now by the Atiyah-Floer conjecture there should be a natural isomorphism between the
topological invariant $\HF^{\rm inst}_*(H_0\cup_\S H_1)$ and the symplectic invariant
$\HF^{\rm symp}_*(\cR_\S,L_{H_0},L_{H_1})$ -- assuming that the first is defined, 
i.e.\ $H_0\cup H_1$ is a homology $3$-sphere, and that the second can be defined 
in spite of the singularities.
On the gauge theoretic side one obtains two smooth (though infinite dimensional)
Lagrangian submanifolds $\cL_{H_0},\cL_{H_1}\subset\cA(\S)$, to which we can
associate the new invariant 
$\HF^{\rm inst}_*([0,1]\times\S,\cL_{H_0}\times\cL_{H_1})$.
This invariant is more generally defined in the setting below, where 
we again fix $\rG=\SU(2)$.
Here we replace $[0,1]\times\S$ by a more general $3$-manifold $Y$ with boundary
with boundary $\pd Y=\S$.
Then for a union of handle bodies $H=\bigsqcup H_i$
with boundary $\pd H=\bigsqcup \S_i=\S$
we denote by $\cL_H\subset\cA(\S)$ the Lagrangian submanifold
${\cL_{H_0}\times\dots\times\cL_{H_N}\subset\cA(\S_0)\times\dots\times\cA(\S_N)}$.

\begin{thm} {\bf (\cite{SW}) } \label{thm:HF}
Let $Y$ be a compact, oriented $3$-manifold with boundary $\S$.
Let $H$ be a disjoint union of handle bodies with $\pd H=\S$,
and suppose that $Y\cup_\S H$ is a homology $3$-sphere (with $\Z$-coefficients).
Then the Floer homology $\HF^{\rm inst}_*(Y,\cL_H)$ is well-defined
and independent of the metric and perturbations of 
(\ref{new crit}) and (\ref{new flow}) used to define it.
\end{thm}

In this setting, Floer's original invariant $\HF^{\rm inst}_*(Y\cup_\S H)$ 
is also defined, and we expect our invariant to carry the same information.

\begin{con} \label{con:inst}
There is a natural isomorphism
$$
\HF^{\rm inst}_*(Y,\cL_H) \cong \HF^{\rm inst}_*(Y\cup_\S H) .
$$
\end{con}

Hence the new Floer homology with Lagrangian boundary conditions
fits into the Atiyah-Floer conjecture as well 
as for an approach to defining an invariant for more general $3$-manifolds.
In the next section we explain its definition in more detail for the
model case $Y=[0,1]\times\S$ and ${\cL_H=\cL_{H_0}\times\cL_{H_1}}$,
which also is the relevant case for the Atiyah-Floer conjecture.

\section{Instanton and symplectic Floer homologies}
\label{sec:FH}

This section sketches the instanton and symplectic versions of Floer theory 
and compares the analytic behaviour of the underlying trajectory equations.
The purpose of this is to explain the definition of the new 
instanton Floer homology with Lagrangian boundary conditions (in-L)
and to show how it fits between the instanton Floer homology (inst) and the
symplectic Floer homology (symp) and thus provides an intermediate
invariant for approaching the Atiyah-Floer conjecture.
In fact, its trajectories exhibit this interpolation between
anti-self-dual instantons (in their interior behaviour) 
and pseudoholomorphic curves (in their semiglobal behaviour at the boundary).

\smallskip

\noindent
{\bf (inst):}
Let $Y$ be a homology $3$-sphere, i.e.\ a compact oriented
$3$-manifold with integer homology $H_*(Y,\Z)\cong H_*(S^3,\Z)$.
The instanton Floer homology $\HF^{\rm inst}_*(Y)$ was defined by Floer~\cite{F1}.
The basic analytic results for this setup that will be quoted below are mainly
due to Uhlenbeck~\cite{U1,U2}.

\smallskip

\noindent
{\bf (in-L):}
Let $Y=H_0\cup_\S H_1$ be the Heegard splitting of a homology $3$-sphere
into two handle bodies $H_0,H_1$ with common boundary $\pd H_i=\S$. 
We describe the special case 
$\HF^{\rm inst}_*([0,1]\times\S,\cL_{H_0}\times\cL_{H_1})$
of the new instanton Floer homology with Lagrangian boundary conditions 
of theorem~\ref{thm:HF}.
The analytic results for this case are established in~\cite{W elliptic, W bubbling}.

\smallskip

\noindent
{\bf (symp):}
Let $(M,\o)$ be a compact symplectic manifold and assume that it is 
simply connected, positive ($c_1(\rT M)=\l[\o]$ with $\l>0$), 
and has minimal Chern number $N\geq 2$
(where $\la c_1,\pi_2(M)\ra=N\Z$). 
Let $L_0, L_1 \subset M$ be two simply connected Lagrangian submanifolds.
Then the symplectic Lagrangian intersection Floer homology 
$\HF^{\rm symp}_*(M,L_0,L_1)$ is defined by \cite{F2} and many other authors.
The underlying analytic fact here is Gromov's compactness for pseudoholomorphic
curves~\cite{G}.

\smallskip

The instanton cases use the trivial $\SU(2)$-bundle as before.
In the third case one should think of $M=\cR_\S$ and $L_i=L_{H_i}$.
However, their Floer homology is not yet well-defined
due to the quotient singularities.
We do not give complete definitions of the Floer homologies here. 
More detailed expositions can be found in e.g.~\cite{Donaldson book, Sa lecture}. 
In particular, we do not mention the necessary perturbations of the equations for
critical points and trajectories.

\begin{dfn} \label{def:crit}
A {\bf critical point} is 

\medskip

\noindent
{\bf (inst):}
a flat connection $A\in\cA_{\rm flat}(Y)$.

\smallskip

\noindent
{\bf (in-L):}
a flat connection $A+\Psi\dt\in\cA_{\rm flat}([0,1]\times\S)$
with Lagrangian

\qquad\, boundary conditions $A(j)\in\cL_{H_j}$ for $j=0,1$.

\smallskip

\noindent
{\bf (symp):}
an intersection point $x\in L_0\cap L_1$.
\end{dfn}

\noindent
In all three cases, the Floer chain complex is generated by the critical points,
$$
{\rm CF}_* = \bigoplus_{x\;\text{crit.pt.}} \Z \la x \ra .
$$
(In the two instanton cases the generators actually are gauge equivalence classes 
$x=[A]$ or $x=[A+\Psi\dt]$, and the trivial connection is disregarded.)
The boundary operator $\pd:{\rm CF}_* \to {\rm CF}_*$ is defined by counting
trajectories,
$$
\pd \la x^- \ra = \sum_{x^+\,\text{crit.pt.}} \# \cM^0(x^-,x^+) \; \la x^+ \ra .
$$
Here $\cM^0(x^-,x^+)$ is the $0$-dimensional part of the space of trajectories
from $x^-$ to $x^+$. This will be a smooth, compact, oriented manifold, so its 
points can be counted with signs.
The trajectory equations will be given below for the three cases.
The main issue that we then discuss is the compactness of the space 
of trajectories, which will allow the definition of $\pd$.
To obtain a chain complex, one moreover has to establish $\pd\comp\pd=0$
by identifying the boundary of the $1$-dimensional part of the space of 
trajectories with the broken trajectories that contribute to $\pd\comp\pd$.
The Floer homology in the different cases then is the 
homology $H_*({\rm CF},\pd)$ of the corresponding Floer chain complex.
It is graded modulo $8$ in the instanton cases and modulo $2N$ in the
symplectic case.

The trajectory equation depends on the choice of auxiliary data,
that the Floer homology will not depend on.
In the instanton cases this is a metric on $Y$ or $[0,1]\times\S$ respectively.
(In the second case we will give the equation for a product metric.)
In the symplectic case we fix an $\o$-compatible almost 
complex structure $J$ on $M$.
The moduli space of trajectories then is the space of solutions
of the trajectory equation modulo time shift (in the $\R$-variable) and modulo
gauge equivalence in the instanton cases.

\begin{dfn} \label{def:trajectory}
A {\bf trajectory} is a solution of the trajectory equation $(T)$.

\medskip

\noindent
{\bf (inst):}
An anti-self-dual instanton on $\R\times Y$:
\begin{align*}
&B:\R\to\cA(Y) \quad\text{satisfying} \\
&(T)\qquad \pd_s B + * F_B = 0  \qquad\qquad\qquad\qquad\qquad\qquad\quad
\end{align*}

\smallskip

\noindent
{\bf (in-L):}
An anti-self-dual instanton on $\R\times[0,1]\times\S$
with Lagrangian 

\qquad\, boundary conditions:
\begin{align*} 
& (A,\Psi):\R\times[0,1]\to\cA(\S)\times\cC^\infty(\S,\su(2))
\quad\text{satisfying} \\
&(T)\qquad \left\{ \begin{aligned}
&\pd_s A + * ( \pd_t A - \rd_A\Psi ) =0 \\
&\pd_s \Psi + * F_A  =0 \\
&A(s,j) \in \cL_{H_j} \qquad\forall s\in\R , j\in\{0,1\}
\end{aligned} \right.
\end{align*}

\smallskip

\noindent
{\bf (symp):}
A $J$-holomorphic strip with Lagrangian boundary conditions:
\begin{align*}
&u:\R\times[0,1]\to M \quad\text{satisfying} \\
&(T)\qquad \left\{ \begin{aligned}
&\pd_s u + J \pd_t u =0 \\
&u(s,j) \in L_j \qquad\forall s\in\R , j\in\{0,1\} \qquad
\end{aligned} \right.
\end{align*}
\end{dfn}

Pictures of these trajectories and a table that summarizes the definitions and results 
for the three Floer theories can be found on page~\pageref{pic:trajectories}
and~\pageref{table 2}. 
The equation in case (in-L) is $\pd_s B + * F_B = 0$ for
$B=A+\Psi\dt$, and in both instanton cases this is the anti-self-duality
equation for the connection $\X=0\ds+B$ in temporal gauge; c.f.\ section~\ref{sec:gauge}.

To ensure that the trajectories converge to critical points as the
$\R$-variable tends to $\pm\infty$, one needs some a priori bound.
This is provided by energy functionals given in the lemma below
(a consequence of theorems~\ref{thm:C},~\ref{thm:EQ}).

\begin{lem}
If a trajectory has finite {\bf energy} $\cE$, then it converges (exponentially)
to critical points as $\R\ni s\to\pm\infty$.

\medskip

\noindent
{\bf (inst):}

\vspace{-11mm}

\begin{align*}
&\cE(B)=\int_{\R\times Y} |\pd_s B|^2 <\infty \\
&\Rightarrow  B(s) \underset{s\to\pm\infty}{\longrightarrow} B^\pm \in\cA_{\rm flat}(Y)
\qquad\qquad\qquad\qquad\qquad\quad
\end{align*}

\smallskip

\noindent
{\bf (in-L):}

\vspace{-11mm}

\begin{align*}
&\cE(A,\Psi)=\int_{\R\times[0,1]\times\S} |\pd_s A_i|^2 +|F_{A_i}|^2 <\infty \\
&\Rightarrow  A(s)+\Psi(s)\dt \underset{s\to\pm\infty}{\longrightarrow}
A^\pm + \Psi^\pm\dt \in\cA_{\rm flat}([0,1]\times Y);\\
&\qquad\qquad\qquad\qquad\qquad\quad
A^\pm(0)\in\cL_{H_0}, A^\pm(1)\in\cL_{H_1}
\end{align*}

\smallskip

\noindent
{\bf (symp):} 

\vspace{-11mm}
\begin{align*}
&\cE(u)=\int_{\R\times[0,1]} |\pd_s u |^2 <\infty \\
&\Rightarrow  u(s,\cdot) \underset{s\to\pm\infty}{\longrightarrow} x^\pm \in L_0\cap L_1
\qquad\qquad\qquad\qquad\qquad
\end{align*}
\end{lem}

In the two instanton cases, the energy of a trajectory
equals to the Yang-Mills energy $\half \int |F_\X|^2$ 
of the corresponding anti-self-dual connection.
In all cases the energy is conformally invariant, so by rescaling one solution
one can obtain a sequence of solutions (on a ball) whose energy is bounded,
but that blows up at one point -- where all the energy concentrates.
This effect can be excluded by assuming that the energy density does not blow up.
For all three equations, this is enough to obtain $\cC^\infty_{\rm loc}$-compactness.

\begin{thm}{\bf (Compactness)} \label{thm:C}
Consider a sequence of trajectories and
suppose that their {\bf energy density} is locally uniformly bounded:

\smallskip

\noindent
{\bf (inst):}
$|\pd_s B_i|^2$ is locally uniformly bounded on $\R\times Y$.

\smallskip

\noindent
{\bf (in-L):}
$\|\pd_s A_i\|_{L^2(\Sigma)}^2 +\|F_{A_i}\|_{L^2(\Sigma)}^2$ 
is locally uniformly bounded on $\R\times[0,1]$.

\smallskip

\noindent
{\bf (symp):}
$|\pd_s u_i|^2$ is locally uniformly bounded on $\R\times[0,1]$.

\smallskip

\noindent
Then, after going to a subsequence, and in the cases (inst), (in-L) applying 
a sequence of gauge transformations $g_i\in\cG(\R\times Y)$ or 
$g_i\in\cG(\R\times [0,1]\times\S)$,
the trajectories converge uniformly with all derivatives on every compact 
subset (i.e.\ in the $\cC^\infty_{\rm loc}$-topology) to a new trajectory.
\end{thm}

The compactness statement in case (in-L) in fact also holds when the Lagrangians
$\cL_{H_i}$ are replaced with general gauge invariant Lagrangians as on 
page~\pageref{Lag ass}.
This result was proven in \cite{W elliptic} under the (stronger) 
standard assumption from gauge theory that 
${|\pd_s B_i|^2=|\pd_s A_i|^2+|F_{A_i}|^2}$ is locally uniformly 
bounded on $\R\times[0,1]\times\S$ (or is locally $L^p$-bounded for $p>2$).
The weaker assumption above implies pointwise bounds in the interior
by a mean value inequality.
Near the boundary this is not a direct consequence, but an extra argument 
\cite[Lemma~2.4]{W bubbling} 
provides local $L^p$-bounds for any $p<3$.
Thus we can state the compactness result in this form, which already hints
at a similar behaviour to pseudoholomorphic curves on $\R\times[0,1]$.
This stronger statement becomes crucial in the bubbling analysis below.

The goal of our analytic discussion of the trajectory equation is to understand
the compactness or compactification of the $k$-dimensional part $\cM^k(x^-,x^+)$ 
of the space of trajectories with fixed limits $x^\pm$.
(Here $k=0$ and $k=1$ are relevant for the definition of $\pd$ and for the proof of $\pd\comp\pd=0$.)

The assumptions in theorem~\ref{thm:C} are too strong for that purpose since 
we only have a bound on the energy, not on the energy density, of trajectories in $\cM^k(x^-,x^+)$.
In fact, in the three present cases the energy of a trajectory is uniquely determined by 
its limits $x^-,x^+$ and its index $k$ via a monotonicity formula.
So we need to consider a sequence of trajectories with fixed energy
and analyze the possible divergence of the sequence when the uniform bounds 
in theorem~\ref{thm:C} do not hold.
This divergence is usually described by the 'bubbling off' of some part of the trajectory:
In the case (inst) the 'bubbles' are instantons on $S^4$; in the case (symp) they are 
pseudoholomorphic spheres or disks.
In the new case (in-L) we also encounter instantons on $S^4$ 'bubbling off' at both interior
or boundary points. 
Additional 'bubbles' in the form of anti-self-dual instantons on the half space 
were expected in \cite{Sa}. Our result below now seems to indicate a semiglobal 
bubbling effect at the boundary, which conjecturally might be described 
as a holomorphic disk in the space of connections $\cA(\S)$.
Fortunately, the geometric understanding of the bubbles is not necessary 
for the purpose of Floer theory in the monotone case. 
It can be replaced by an analytic understanding of the bubbling
in the form of the following energy quantization result.

For the purpose of this statement we abbreviate $Y=[0,1]\times\S$ in case (in-L)
and $Y=[0,1]$ in case (symp), so all trajectories are defined on $\R\times Y$.

\begin{thm}{\bf (Energy Quantization)} \label{thm:EQ}
There exists a constant $\hbar>0$ such that the following holds.
Consider a sequence of trajectories whose energy is bounded by some $E<\infty$.

Then, after going to a subsequence, the energy densities are locally 
uniformly bounded as in theorem~\ref{thm:C} 
on $(\R\times Y)\setminus\bigcup_{k=1}^N P_k$, 
the complement of a finite union 
of bubbling loci $P_k$ as below.
At each bubbling locus $P_k$ there is a concentration of energy of at 
least $\hbar$ on neighbourhoods with radii $\ep_i\to 0$.

\medskip

\noindent
{\bf (inst):}
Each bubbling locus is a point $P_k=x_k\in\R\times Y$ with 
$$
\int_{B_{\ep_i}(x_k)} |\pd_s B_i|^2 \geq \hbar .
$$

\smallskip

\noindent
{\bf (in-L):}
Each bubbling locus is 

\qquad\, 
either an interior point ${P_k=x_k\in\R\times(0,1)\times\Sigma}$ with
$$
\int_{B_{\ep_i}(x_k)} |\pd_s A_i|^2 + |F_{A_i}|^2 \geq \hbar ,
$$

\qquad\, or a boundary slice $P_k=\{(s_k,t_k)\}\times\Sigma$,
$(s_k,t_k)\in\R\times\{0,1\}$ with
$$
\int_{B_{\ep_i}(s_k,t_k)} \|\pd_s A_i\|_{L^2(\S)}^2 
+ |F_{A_i}|_{L^2(\S)}^2 \geq \hbar .
$$

\smallskip

\noindent
{\bf (symp):} 
Each bubbling locus is a point $P_k=(s_k,t_k)\in\R\times[0,1]$ with
$$
\int_{B_{\ep_i}(s_k,t_k)} |\pd_s u_i|^2 \geq \hbar .
$$
\end{thm}

In case (in-L) both an instanton on $S^4$ bubbling off at a boundary point
and the conjectural holomorphic disk in $\cA(\S)$ are 
described by a boundary slice as bubbling locus.
The proof in case (in-L) goes along the lines of an energy quantization principle
explained in \cite{W mean} but deals with some additional difficulties.
In the cases (inst) and (symp) the above result can be obtained straight forward
from this principle and a control on the Laplacian (and normal derivative) 
of the energy density. See section~\ref{sec:analysis} for details.

The combination of theorems~\ref{thm:C} and~\ref{thm:EQ} can be rephrased as:
'There is a $\cC^\infty_{\rm loc}$-convergent subsequence if the energy is locally small.'
In the cases (inst) and (symp) it is sufficient to assume that every point
in $\R\times Y$ or $\R\times[0,1]$ respectively
has a neighbourhood on which the energy of each trajectory in 
the sequence is less than $\hbar$.
In the case (in-L) this assumption is the same for points in the interior
$\R\times(0,1)\times\S$. For a point $(s,j,z)\in\R\times\{0,1\}\times\S$ 
on the boundary however, it is not enough to assume that the energies are small
on a neighbourhood of that point, but one needs to assume that 
there is a neighbourhood of the whole boundary slice $\{(s,j)\}\times\S$
on which the energy of each trajectory in the sequence is less than $\hbar$.

The full consequence of theorems~\ref{thm:EQ} and~\ref{thm:C}
is the following compactness.

\begin{cor}\label{cor:comp}
Consider a sequence of trajectories with energy bounded by $E<\infty$.
Then, after going to a subsequence, there exist finitely many bubbling loci 
$P_1,\ldots,P_N$ as in theorem~\ref{thm:EQ}, 
and in the cases (inst) and \hbox{(in-L)} there exists a sequence of gauge 
transformations in $\cG((\R\times Y)\setminus\bigcup_{i=1}^k P_k)$, such that 
the trajectories converge (after gauge transformation) in the $\cC^\infty_{\rm loc}$-topology 
on ${(\R\times Y)\setminus\bigcup_{i=1}^k P_k}$ to a new solution of the trajectory
equation (T) on ${(\R\times Y)\setminus\bigcup_{i=1}^k P_k}$ with energy 
$\cE\leq E - N\hbar$.
\end{cor}

Keep in mind that the bubbling loci $P_k$ and thus the singularities of the new
solution obtained in corollary~\ref{cor:comp} are always points, except for the
case (in-L) where $2$-dimensional singularities can occur at the boundary.
The next step in the compactification (or proof of compactness) of the spaces 
of trajectories is to remove these singularities.
We give a general statement that is a consequence of the subsequent
removable singularity theorems for the local models of the singularities.

Here $B^n$ denotes the unit ball in $\R^n$ centered at $0$, 
and ${D^2:=B^2\cap\H^2}$ is the unit half ball in  
the half space ${\H^2=\{(s,t)\in\R^2 \st t\geq 0 \}}$ with center~$0$.
In~the two boundary cases, the Lagrangian submanifold $\cL_H$ or $L$ can be
either of the two $\cL_{H_i}$ or $L_i$ respectively.

\begin{thm}{\bf (Removal of Singularities)} \label{thm:RemSing}
Consider a smooth solution of the trajectory equation (T) on 
$(\R\times Y)\setminus \bigcup_{k=1}^N P_k$ that has finite 
energy~$E$.
Then (in case (inst) and (in-L) after applying a gauge transformation
in $\cG((\R\times Y)\setminus \bigcup_{k=1}^N P_k)$) 
the solution extends to a trajectory on $\R\times Y$ with energy $E$.

\medskip
\noindent
{\bf (inst), (in-L,interior):}
Suppose that $\X\in\cA(B^4\setminus\{0\})$ satisfies
$$
\qquad F_\X+*F_\X=0
\qquad\text{and}\qquad
\int_{B^4\setminus\{0\}} |F_\X|^2 <\infty .
$$
Then there exists a gauge transformation $g\in\cG(B^4\setminus\{0\})$ such that
$g^*\X$ extends to a solution $\tilde\X\in\cA(B^4)$.

\medskip
\noindent
{\bf (in-L,boundary):}
Suppose that $\X\in\cA((D^2\setminus\{0\})\times\S)$ satisfies
\[
\qquad \left\{\begin{aligned}
&F_\X+*F_\X=0 \\
&\X|_{\{(s,0)\}\times\S}\in\cL_{H} \quad\forall s 
\end{aligned}\right.
\qquad\text{and}\qquad
\int_{D^2\setminus\{0\}} \int_\S |F_\X|^2 <\infty .
\]
Then there exists a gauge transformation 
$g\in\cG((D^2\setminus\{0\})\times\S)$ such that
$g^*\X$ extends to a solution $\tilde\X\in\cA(D^2\times\S)$.

\medskip
\noindent
{\bf (symp,boundary):}
Suppose that $u\in\cC^\infty(D^2\setminus\{0\},M)$ satisfies
\[
\left\{\begin{aligned}
& \pd_s u + J\pd_t u=0 \\
& u(s,0)\in L \quad\forall s 
\end{aligned}\right.
\qquad\text{and}\qquad
\int_{D^2\setminus\{0\}} |\pd_s u|^2 <\infty .
\]
Then $u$ extends to a solution $\tilde u\in\cC^\infty(D^2,M)$.

\medskip
\noindent
{\bf (symp,interior):}
Suppose that $u\in\cC^\infty(B^2\setminus\{0\},M)$ satisfies 
$$
\pd_s u + J\pd_t u=0
\qquad\text{and}\qquad
\int_{B^2\setminus\{0\}} |\pd_s u|^2 <\infty .
$$
Then $u$ extends to a solution $\tilde u\in\cC^\infty(B^2,M)$.
\end{thm}

In the case (in-L) Uhlenbeck's removable singularity theorem~\cite{U1}
applies to the bubbling loci in the interior. At the boundary
we have to remove $2$-dimensional singularities of an anti-self-dual instanton.
In the interior there would be an obstruction to removing such singularities: 
The holonomies of small loops around the singularity might have a nontrivial limit.
%
%
So it is important to note that this 'pseudoholomorphic behaviour' of the (in-L)
trajectories only occurs at the boundary, where one does not have an obstruction
since there are no loops around the singularity.
One can then imitate the removal of the singularity of a 
pseudoholomorphic curve on $D^2\setminus\{0\}$ with Lagrangian boundary conditions 
to remove the singularity of an anti-self-dual instanton on $(D^2\times\{0\})\times\S$.
This uses an isoperimetric inequality for a local Chern-Simons functional 
instead of the local symplectic action. So far, the definition of this local 
Chern-Simons functional crucially uses the fact that the Lagrangian boundary 
condition arises from a handle~body.

The final result of the analysis of trajectories in theorems~\ref{thm:C},
~\ref{thm:EQ}, and~\ref{thm:RemSing} is that the moduli spaces of trajectories
are compact up to 'bubbling' and 'breaking of trajectories'.
Here 'bubbling' means the concentration of energy at a bubbling locus 
as in theorem~\ref{thm:EQ}.
The 'breaking of trajectories' occurs when a sequence of trajectories with 
constant energy converges smoothly on every compact set to a new trajectory, 
but the limit has less energy. In that case, the energy difference must have
moved out to $s\to\pm\infty$ and can be recaptured as the energy of a limit of
shifted trajectories.
A standard iteration of such shifts yields a finite collection of trajectories 
(a 'broken trajectory') whose total energy equals to the fixed energy 
of the sequence.

To proceed with the definition of $\pd$ and the proof of $\pd\comp\pd=0$ 
one needs to perturb the trajectory equation (T) so that the moduli spaces 
$\cM^k(x^-,x^+)$ of trajectories become smooth manifolds. 
Here a priori $k\in\Z$ is the index of a Fredholm operator (the linearization of (T)) 
associated to the trajectories. 
For a smooth moduli space, $k$ equals to the dimension of the component,
hence $\cM^k(x^-,x^+)$ is empty for $k\leq-1$.
By a monotonicity formula, $k$ moreover determines the energy of the trajectories
such that a trajectory of lower energy has to lie in a moduli space of lower dimension.
From this one can deduce that $\cM^0(x^-,x^+)$ is compact (and thus can be counted to
define~$\pd$): 
It consists of trajectories with the minimal energy that allows
to connect $x^-$ to $x^+$. So bubbling can be ruled out since (after removal of the
singularities) it would lead to a trajectory of even lower energy.
The breaking of trajectories is ruled out by a similar index-energy argument.

Bubbling is also excluded in $\cM^k(x^-,x^+)$ for $k\leq 7$ 
(or $2N-1$ in the symplectic case) since $x^-$ and 
$x^+$ determine the index $k$ modulo $8$ (or $2N$). So a loss of energy
corresponds to a jump by $8$ (or $2N$) in the dimension.
The breaking of trajectories is no longer ruled out; on the contrary,
${\pd\comp\pd=0}$ follows from the fact that the ends of the $1$-dimensional moduli 
spaces exactly correspond to the broken trajectories which are counted by $\pd\comp\pd$.

\begin{landscape}
\small
\begin{tabular}[b]{p{0.4\textwidth}|p{0.5\textwidth}|p{0.45\textwidth}}
{\begin{minipage}[t]{0.4\textwidth}
\label{table 1}
$$\text{\underline{3-manifold topology}}$$
\end{minipage} } 
&
{\begin{minipage}[t]{0.45\textwidth}
$$\text{\underline{gauge theory}}$$
\end{minipage} } & 
{\begin{minipage}[t]{0.5\textwidth}
$$\text{\underline{symplectic topology}}$$
\end{minipage} } 
\\
{\begin{minipage}[t]{0.4\textwidth}
\begin{center}
\vspace{3mm}
$\S$ \\
Riemann surface
\end{center}
\end{minipage} } 
&
{\begin{minipage}[t]{0.45\textwidth}
\begin{center}
\vspace{3mm}
$\cA(\S)$  \\
symplectic Banach space
\end{center}
\end{minipage} } 
&
{\begin{minipage}[t]{0.5\textwidth}
\begin{align*}
&\cR_\S = \cA_{\rm flat}(\S)/\cG(\S) = \cA(\S)/\hspace{-1mm}/\cG(\S)  \\
& \qquad \cong \Hom(\pi_1(\S),\SU(2))/\SU(2) \\
&\text{(singular) symplectic manifold}
\end{align*}
\end{minipage} } 
\\
{\begin{minipage}[t]{0.4\textwidth}
\begin{center}
\vspace{3mm}
$H , \quad \pd H = \S$ \\
handle body
\end{center}
\end{minipage} }
&
{\begin{minipage}[t]{0.45\textwidth}
\begin{center}
\vspace{3mm}
$\cL_H = \cA_{\rm flat}(H)|_\S$ \\
Lagrangian Banach submanifold
\end{center}
\end{minipage} } 
& 
{\begin{minipage}[t]{0.5\textwidth}
\begin{align*}
&L_H  = \cA_{\rm flat}(H)|_\S/\cG(\S) \\
&\cong \Hom(\pi_1(\S)/\pd\pi_2(H,\S),\SU(2))/\SU(2) \\
&\text{(singular) Lagrangian submanifold}
\end{align*}
\end{minipage} }
\\
{\begin{minipage}[t]{0.4\textwidth}
\begin{center}
\vspace{3mm}
$Y=H_0\cup_\S H_1$ \\
Heegard splitting
\end{center}
\end{minipage} } 
&
{\begin{minipage}[t]{0.45\textwidth}
\begin{align*}
&\cL_{H_0},\cL_{H_1} \subset \cA(\S)
\end{align*}
\end{minipage} } 
& 
{\begin{minipage}[t]{0.5\textwidth}
\begin{align*}
&L_{H_0},L_{H_1}\subset\cR_\S
\end{align*}
\end{minipage} } 
\\
{\begin{minipage}[t]{0.4\textwidth}
\begin{align*}
& \HF^{\rm inst}_*(H_0\cup_\S H_1)
\end{align*}
\end{minipage} } 
&
{\begin{minipage}[t]{0.45\textwidth}
\begin{align*}
&\HF^{\rm inst}_*([0,1]\times\S,\cL_{H_0}\times\cL_{H_1})
\end{align*}
\end{minipage} } 
& 
{\begin{minipage}[t]{0.5\textwidth}
\begin{align*}
&\HF^{\rm symp}_*(\cR_\S,L_{H_0},L_{H_1})
\end{align*}
\end{minipage} } 
\\
\begin{center}
\begin{picture}(0,0)%
\includegraphics{trajectory1.pstex}%
\end{picture}%
\setlength{\unitlength}{2368sp}%
\begingroup\makeatletter\ifx\SetFigFont\undefined%
\gdef\SetFigFont#1#2#3#4#5{%
  \reset@font\fontsize{#1}{#2pt}%
  \fontfamily{#3}\fontseries{#4}\fontshape{#5}%
  \selectfont}%
\fi\endgroup%
\begin{picture}(2769,2752)(924,-2633)
\put(3076,-1261){\makebox(0,0)[lb]{\smash{\SetFigFont{10}{12.0}{\familydefault}{\mddefault}{\updefault}$H_1$}}}
\put(2401,-1336){\makebox(0,0)[lb]{\smash{\SetFigFont{10}{12.0}{\familydefault}{\mddefault}{\updefault}$\Sigma$}}}
\put(2401,-2386){\makebox(0,0)[lb]{\smash{\SetFigFont{10}{12.0}{\familydefault}{\mddefault}{\updefault}$Y$}}}
\put(1576,-1261){\makebox(0,0)[lb]{\smash{\SetFigFont{10}{12.0}{\familydefault}{\mddefault}{\updefault}$H_0$}}}
\put(1126,-61){\makebox(0,0)[lb]{\smash{\SetFigFont{10}{12.0}{\familydefault}{\mddefault}{\updefault}$\mathbb{R}$}}}
\end{picture}

\end{center}
\label{pic:trajectories}
&
\begin{center}
\begin{picture}(0,0)%
\includegraphics{trajectory2.pstex}%
\end{picture}%
\setlength{\unitlength}{2368sp}%
\begingroup\makeatletter\ifx\SetFigFont\undefined%
\gdef\SetFigFont#1#2#3#4#5{%
  \reset@font\fontsize{#1}{#2pt}%
  \fontfamily{#3}\fontseries{#4}\fontshape{#5}%
  \selectfont}%
\fi\endgroup%
\begin{picture}(3421,2752)(1224,-5708)
\put(2626,-5461){\makebox(0,0)[lb]{\smash{\SetFigFont{10}{12.0}{\familydefault}{\mddefault}{\updefault}$[0,1]\times\Sigma$}}}
\put(1426,-3136){\makebox(0,0)[lb]{\smash{\SetFigFont{10}{12.0}{\familydefault}{\mddefault}{\updefault}$\mathbb{R}$}}}
\put(1726,-4336){\makebox(0,0)[lb]{\smash{\SetFigFont{10}{12.0}{\familydefault}{\mddefault}{\updefault}$\mathcal{L}_{H_0}$}}}
\put(4051,-4336){\makebox(0,0)[lb]{\smash{\SetFigFont{10}{12.0}{\familydefault}{\mddefault}{\updefault}$\mathcal{L}_{H_1}$}}}
\end{picture}

\end{center}
&
\begin{center}
\begin{picture}(0,0)%
\includegraphics{trajectory3.pstex}%
\end{picture}%
\setlength{\unitlength}{2368sp}%
\begingroup\makeatletter\ifx\SetFigFont\undefined%
\gdef\SetFigFont#1#2#3#4#5{%
  \reset@font\fontsize{#1}{#2pt}%
  \fontfamily{#3}\fontseries{#4}\fontshape{#5}%
  \selectfont}%
\fi\endgroup%
\begin{picture}(3462,2649)(151,-2623)
\put(151,-661){\makebox(0,0)[lb]{\smash{\SetFigFont{10}{12.0}{\familydefault}{\mddefault}{\updefault}$\mathcal{R}_\Sigma$}}}
\put(2626,-1261){\makebox(0,0)[lb]{\smash{\SetFigFont{10}{12.0}{\familydefault}{\mddefault}{\updefault}$L_{H_1}$}}}
\put(1051,-1261){\makebox(0,0)[lb]{\smash{\SetFigFont{10}{12.0}{\familydefault}{\mddefault}{\updefault}$L_{H_0}$}}}
\end{picture}

\end{center}
\end{tabular}
\end{landscape}

\begin{landscape}
\small\label{table 2}
\begin{tabular}[b]{p{0.21\textwidth}|p{0.25\textwidth}|p{0.5\textwidth}|p{0.5\textwidth}}
& 
{\begin{minipage}[t]{0.25\textwidth}
$$\underline{\HF^{\rm inst}_*(Y)}$$
\end{minipage} } 
& 
{\begin{minipage}[t]{0.5\textwidth}
$$\underline{\HF^{\rm inst}_*([0,1]\times\S,\cL_{H_0}\times\cL_{H_1})}$$ 
\end{minipage} } 
&
{\begin{minipage}[t]{0.45\textwidth}
$$\underline{\HF^{\rm symp}_*(M,L_0,L_1)}$$ 
\end{minipage} } \\
{\begin{minipage}[t]{0.21\textwidth}
\vspace{3.5mm}
'critical points' :
\end{minipage} } 
&
{\begin{minipage}[t]{0.25\textwidth}
\begin{align*}
B\in\cA_{\rm flat}(Y)
\end{align*}
\end{minipage} }
& 
{\begin{minipage}[t]{0.5\textwidth}
\begin{align*}
A+\Psi\dt  &\in\cA_{\rm flat}([0,1]\times\S) ,\\ 
A(j)|_{\pd Y} &\in\cL_{H_j}
\end{align*}
\end{minipage} }
&
{\begin{minipage}[t]{0.45\textwidth}
\begin{align*}
x\in L_0\cap L_1
\end{align*}
\end{minipage} }
\\
{\begin{minipage}[t]{0.21\textwidth}
\vspace{3.5mm}
trajectories :
\end{minipage} }
&
{\begin{minipage}[t]{0.25\textwidth}
\begin{align*} 
& B:\R\to\cA(Y) \\
& \Bigl\{ \pd_s B + * F_B = 0 \\
\end{align*}
\end{minipage} }
&
{\begin{minipage}[t]{0.5\textwidth}
\begin{align*} 
& (A,\Psi):\R\times[0,1]\to\cA(\S)\times\cC^\infty(\S,\su(2))  \\
&\left\{ \begin{aligned}
&\pd_s A + * ( \pd_t A - \rd_A\Psi ) =0 \\
&\pd_s \Psi + * F_A  =0 \\
&A(s,j) \in \cL_{H_j} 
\end{aligned} \right.
\end{align*}
\end{minipage} }
&
{\begin{minipage}[t]{0.45\textwidth}
\begin{align*} 
& u:\R\times[0,1]\to M \\
&\left\{ \begin{aligned}
& \pd_s u + J \pd_t u =0 \\
& u(s,j) \in L_j 
\end{aligned} \right.
\end{align*}
\end{minipage} }
\\
{\begin{minipage}[t]{0.21\textwidth}
\vspace{5mm}
energy :
\end{minipage} }
&
{\begin{minipage}[t]{0.25\textwidth}
$$\int_{\R\times Y} |\pd_s B|^2 \leq C$$
\end{minipage} }
&
{\begin{minipage}[t]{0.5\textwidth}
$$\int_{\R\times [0,1]\times\S} |\pd_s A|^2 + |F_A|^2 \leq C$$
\end{minipage} }
&
{\begin{minipage}[t]{0.45\textwidth}
$$\int_{\R\times [0,1]} |\pd_s u|^2 \leq C$$
\end{minipage} }
\\
{\begin{minipage}[t]{0.21\textwidth}
\vspace{2mm}
uniform bounds for compactness~:
\end{minipage} }
&
{\begin{minipage}[t]{0.25\textwidth}
$$
\sup_{\R\times Y} |\pd_s B|^2 < \infty
$$
\end{minipage} }
&
{\begin{minipage}[t]{0.5\textwidth}
$$
\sup_{\R\times [0,1]} \|\pd_s A\|_{L^2(\S)}^2 + \|F_A\|_{L^2(\S)}^2 <\infty
$$
\end{minipage} }
&
{\begin{minipage}[t]{0.45\textwidth}
$$
\sup_{\R\times [0,1]} |\pd_s u|^2 <\infty
$$
\end{minipage} }
\\
{\begin{minipage}[t]{0.21\textwidth}
\vspace{3.5mm}
bubbling loci :
\end{minipage} }
&
{\begin{minipage}[t]{0.25\textwidth}
\begin{align*}
\text{points}\;\; & x\in\R\times Y
\end{align*}
\end{minipage} }
&
{\begin{minipage}[t]{0.5\textwidth}
\begin{align*}
&\text{interior points}\;\; x\in\R\times (0,1)\times\S, \\
&\text{boundary slices}\;\; \{(s,j)\}\times\S
\end{align*}
\end{minipage} }
&
{\begin{minipage}[t]{0.45\textwidth}
\begin{align*}
\text{interior points}\;\;& (s,t)\in\R\times (0,1), \\
\text{boundary points}\;\;& (s,j)\in\R\times\{0,1\}
\end{align*}
\end{minipage} }
\\
{\begin{minipage}[t]{0.21\textwidth}
\vspace{3mm}
removable \\
singularities :
\end{minipage} }
&
{\begin{minipage}[t]{0.25\textwidth}
\begin{align*}
B^4 \setminus \{0\}
\end{align*}
\end{minipage} }
&
{\begin{minipage}[t]{0.5\textwidth}
\begin{align*}
 B^4 &\setminus \{0\} \\
 (D^2\times\S) &\setminus (\{0\}\times\S)
\end{align*}
\end{minipage} }
&
{\begin{minipage}[t]{0.45\textwidth}
\begin{align*}
& B^2\setminus\{0\} \\
& D^2\setminus\{0\} 
\end{align*}
\end{minipage} }
\end{tabular}
\end{landscape}

\section{The Atiyah-Floer conjecture}
\label{sec:AF}

To give a precise statement of the Atiyah-Floer conjecture we need to
refine the notion of handle bodies and Heegard splittings.
A {\bf handle body} is an oriented $3$-manifold with boundary
that is obtained by attaching finitely many \hbox{$1$-handles} to a $3$-ball.
The {\bf spine} of a handle body $H$ is a graph $S\subset H$ embedded in its 
interior that arises from replacing the ball by a vertex and the handles 
by edges with ends on this vertex.
Its significance is that $H\setminus S\cong [0,1)\times\pd H$, 
so $H$ retracts onto $S$.
For each genus $g\in\N_0$ we fix a standard handle body and
spine $S\subset H$.

\begin{dfn}
A {\bf Heegard splitting} of a closed oriented $3$-manifold $Y$ consists of 
two embeddings $\psi_i:H\hookrightarrow Y$ of a standard handle body $H$
such that $\psi_0(\pd H)=\psi_1(\pd H)=\im\psi_0\cap\im\psi_1$.
We abbreviate the Heegard splitting by ${Y=H_0\cup_\S H_1}$,
where $H_i:=\psi_i(H)\subset Y$ and $\S:=\psi_i(\pd H)=H_1\cap H_2$.
\end{dfn}

Next, a homology $3$-sphere is a compact oriented $3$-manifold $Y$ 
whose integer homology is that of a $3$-sphere, $H_*(Y,\Z)\cong H_*(S^3,\Z)$.

\begin{con} \label{con:AF} {(\bf Atiyah--Floer)}
Let $Y$ be a homology $3$-sphere.
Then every Heegard splitting $Y=H_0\cup_\S H_1$ induces a natural
isomorphism
$$
\HF^{\rm inst}_*(Y) \cong \HF^{\rm symp}_*(\cR_\S,L_{H_0},L_{H_1}).
$$
\end{con}

Here 'natural' in particular means that the isomorphism should be invariant under 
isotopies of the Heegard splitting.
Note that for nonisotopic Heegard splittings of the same genus one can identify the 
$\cR_\S$, but the pairs of Lagrangians (and thus the conjectured
isomorphism) will be different. 
The conjecture would 
then provide isomorphisms between the symplectic Floer homologies arising from 
different Heegard diagrams of the same $3$-manifold.

The first task posed by this conjecture is to give a precise definition of the
symplectic Floer homology for the Lagrangians $L_{H_0},L_{H_1}$ in the singular 
symplectic space $\cR_\S$. They can be viewed as symplectic quotients of the gauge action 
on the smooth Banach-manifolds $\cL_{H_0},\cL_{H_1}\subset\cA(\S)$
(see \cite{AB} and section~\ref{sec:gauge}).
For finite dimensional Hamiltonian group actions, Salamon et al.\ introduced 
invariants based on the symplectic vortex equations on the total space, 
see e.g.\ \cite{CGMS}.
Gaio and Salamon \cite{GS} identified these with the Gromov-Witten 
invariants for smooth and monotone symplectic quotients.
In view of this result, a plausible definition of
$\HF^{\rm symp}_*(\cR_\Si,L_{H_0},L_{H_1})$ could be to replace its ill-defined
trajectories (pseudoholomorophic curves in the singular symplectic quotient)
by solutions of the corresponding symplectic vortex equations:
A triple of maps $A:\R\times[0,1]\to\cA(\S)$ and 
$\P,\Psi:\R\times[0,1]\to\cC^\infty(\S,\su(2))\cong\rT_\smone\cG(\S)$ that satisfy
\begin{equation} \label{eq 1} 
\left\{ \begin{aligned}
( \pd_s A - \rd_A\P ) + * ( \pd_t A - \rd_A\Psi ) &= 0 ,\\
\pd_s \Psi - \pd_t\P + [\P,\Psi] + * F_A  &=0 , \\
 A(s,i) \in \cL_{H_i} \qquad \forall s\in\R, i&\in\{0,1\} .
\end{aligned} \right.
\end{equation}
Here $\P\mapsto\rd_A\P$ is the infinitesimal action and
$A\mapsto *F_A$ is the moment map of the gauge action,
where $*$ is the Hodge operator of a metric $g_\S$ on~$\S$.
This system is the anti-self-duality equation 
with Lagrangian boundary conditions for the connection $\P\ds+\Psi\dt+A$ 
on ${\R\times[0,1]\times\S}$ with respect to the metric $\ds^2+\dt^2+g_\S$,
i.e.\ the trajectory equation of definition~\ref{def:trajectory}
in temporal gauge $\P=0$.
So in this case the symplectic vortex equations lead directly to the 
new Floer homology ${\HF^{\rm inst}_*([0,1]\times\Si,\cL_{H_0}\times \cL_{H_1})}$,
which is well-defined since $([0,1]\times\S) \cup (H_0\sqcup H_1)\cong H_0\cup_\S H_1$ 
is a homology $3$-sphere.
Defining $\HF^{\rm symp}_*(\cR_\Si,L_{H_0},L_{H_1})$ via (\ref{eq 1})
would reduce the Atiyah-Floer conjecture~\ref{con:AF}
to the subsequent special case of conjecture~\ref{con:inst}.
We intend however to give a less far fetched definition of the
symplectic Floer homology and use the following only as first step
towards a proof of the Atiyah-Floer conjecture.

\begin{con}\label{con:1}
Every Heegard splitting $Y=H_0\cup_\Sigma H_1$ of a homology 3-sphere 
induces a natural isomorphism
$$
\HF^{\rm inst}_*(Y)
\cong \HF^{\rm inst}_*([0,1]\times\Si,\cL_{H_0}\times \cL_{H_1}).
$$
\end{con}

To prove this, one has to identify the critical points and trajectories of
both Floer homologies. 
Our idea for a proof uses the following decomposition of $Y$.
We restrict the embeddings $\psi_i$ to the complement of the spine 
${H\setminus S\cong [\half,1)\times\S}$ and glue them at $\{\half\}\times\S$ to obtain
an embedding ${\psi:(0,1)\times\S\hookrightarrow Y}$ 
such that ${\psi(\half,\cdot)={\rm id}_\S}$ and $\psi(t,\S)$ converges 
to the spine $\psi_i(S)\subset H_i$ as $t\to i$ for ${i=0,1}$.
Then 
$$
Y= H_0^\d \sqcup Y_\d \sqcup H_1^\d ;\qquad Y_\d := \psi([\d,1-\d]\times\S) .
$$
Here the $H_i^\d\subset Y$ are isotopic to the open handle bodies ${\rm int}(H_i)$
and ${Y_\d\cong[0,1]\times\S}$ via $\psi\comp\t_\d$, where
$\t_\d:[0,1]\times\S\to[\d,1-\d]\times\S$
is the obvious linear isomorphism.
With this the critical points can be identified elementary as follows:
Every $\tA\in\cA_{\rm flat}(Y)$ can be decomposed and pulled back
to a triple $(A,\tA_0,\tA_1)$ of ${A\in\cA_{\rm flat}([0,1]\times\S)}$ and 
$\tA_i\in\cA_{\rm flat}(H_i)$ such that $A|_{\{i\}\times\S}=\tA_i|_{\pd H_i}$.
So every critical point $[\tA]\in\cR_Y$ corresponds to the 
gauge equivalence class of a flat connection on $[0,1]\times\S$ 
with boundary values in $\cL_{H_0}$ and $\cL_{H_1}$. 
One can check that this in fact gives a bijection between the critical points.
In order to prove conjecture~\ref{con:1} one needs to show that the 
induced map between the Floer complexes is a chain isomorphism.

For that purpose we fix a metric on $Y$ and for a corresponding metric on 
$[0,1]\times\S$ try to establish a bijection between 
the trajectories that contribute to the differential on the two Floer complexes.
(Of course, we have to prove later that the isomorphism is independent of the choices).
A fixed metric on $Y$ gives rise to a family of metrics $g_\d$ on $[0,1]\times\S$
via pullback by ${\psi\comp\t_d:[0,1]\times\S \to Y_\d \subset Y}$.
The metrics $g_\d$ degenerate on $\{0\}\times\S$ and $\{1\}\times\S$ for $\d\to 0$,
but for sufficiently small $\d>0$ we expect to find a bijection between the 
trajectories of $\HF_*^{\rm inst}(Y)$ and those of
${\HF_*^{\rm inst}([0,1]\times\S,\cL_{H_0}\times\cL_{H_1})}$ with respect to~$g_\d$.

The first are anti-self-dual instantons (in temporal gauge) on $\R\times Y$,
that is $B:\R\to\cA(Y)$ satisfying
\begin{align*} 
\pd_s B + * F_B = 0  \qquad\text{on}\;\,\R\times Y .
\end{align*}
The latter are anti-self-dual instantons (in temporal gauge) $A+\Psi\dt$ 
on $\R\times[0,1]\times\S$ with Lagrangian boundary conditions. 
Here the metric $g_\d$ on $[0,1]\times\S$ is not of product form, so the equation (T)
in definition~\ref{def:trajectory} has to be adjusted:
The pair $(A,\Psi)$ is a trajectory if $A+\Psi\dt=\t_\d^*\psi^*B$, where
$B:\R\to\cA(Y_\d)$ is anti-self-dual with respect to
the fixed metric on $Y$ and has boundary values in $\cL_{H_0}$ and $\cL_{H_1}$,
that is
\begin{align*} 
\left\{\begin{aligned}
 \pd_s B + * F_B &= 0  \;\;\,\qquad\text{on}\;\,\R\times Y_\d , \\
B|_{\psi(\{\d\}\times\S)} &= \tB_0  \qquad\text{for some}\;\,\tB_0:\R\to\Af(H_0^\d) , \\
B|_{\psi(\{1-\d\}\times\S)} &= \tB_1 \qquad\text{for some}\;\,\tB_1:\R\to\Af(H_1^\d) .
\end{aligned}\right.
\end{align*}
The task in identifying the trajectories is to consider anti-self-dual 
instantons on $Y_\d$ and transfer between extensions $\tB_i:\R\to\cA(H_i^\d)$
that are slicewise flat ($F_{\tB_i}=0$)
and extensions that are anti-self-dual (${\pd_s \tB_i + *F_{\tB_i}=0}$).
Here the handle bodies $H_i^\d\subset Y$ are small tubes 
around their spines ${\psi_i(S)\subset Y}$.
The restriction of given (anti-self-dual) connections $\tB_i$ to the spines
is up to gauge equivalence determined by their holonomies, 
i.e.\ $SU(2)$-representations of $\pi_1(H_i)$. One can then pick flat connections
on the $H_i^\d$ that have the same holonomy and are close to the $\tB_i$ 
(compared to their energy).
For the converse we will have to use special flat extensions $\tB_i$
with a control on $\pd_s\tB_i$ as in lemma~\ref{lem:ext}.
Combined with the small volume of $H_i^\d$ this should make $\tB_i$ close to anti-self-dual.\\

\begin{center}
\begin{picture}(0,0)%
\includegraphics{degeneration.pstex}%
\end{picture}%
\setlength{\unitlength}{2368sp}%
\begingroup\makeatletter\ifx\SetFigFont\undefined%
\gdef\SetFigFont#1#2#3#4#5{%
  \reset@font\fontsize{#1}{#2pt}%
  \fontfamily{#3}\fontseries{#4}\fontshape{#5}%
  \selectfont}%
\fi\endgroup%
\begin{picture}(9762,2471)(2171,-2458)
\put(10201,-1111){\makebox(0,0)[lb]{\smash{\SetFigFont{10}{12.0}{\familydefault}{\mddefault}{\updefault}${\rm d}s^2 + \epsilon^2 g_\Sigma$}}}
\put(4201,-961){\makebox(0,0)[lb]{\smash{\SetFigFont{10}{12.0}{\familydefault}{\mddefault}{\updefault}$g_\delta$}}}
\put(7201,-1036){\makebox(0,0)[lb]{\smash{\SetFigFont{10}{12.0}{\familydefault}{\mddefault}{\updefault}${\rm d}s^2 + g_\Sigma$}}}
\put(5176,-2386){\makebox(0,0)[lb]{\smash{\SetFigFont{10}{12.0}{\familydefault}{\mddefault}{\updefault}Degenerations of the metric on $[0,1]\times\Sigma$}}}
\end{picture}

\end{center}

The key to this plan of proof is the fact that one can degenerate the metric on 
$[0,1]\times\S$ (as sketched on the left in the above figure) 
without changing the invariant
$\HF^{\rm inst}_*([0,1]\times\Si,\cL_{H_0}\times \cL_{H_1})$.
In the limit of the degeneration one should obtain the invariant 
$\HF^{\rm inst}_*(Y)$ for the closed manifold.
The basic idea of the second step for the Atiyah-Floer conjecture
is to use a second degeneration (on the right in the above sketch)
to transfer from anti-self-dual instantons to pseudoholomorphic 
curves.
This idea was successfully employed by Dostoglou and Salamon \cite{DS}
in their proof of a mapping torus analogon of the Atiyah-Floer conjecture.

A trajectory of the symplectic Floer homology should be a pseudoholomorphic map
$u:\R\times[0,1]\to \cR_\S$ with boundary values in $L_{H_0}$ and~$L_{H_1}$,
\begin{equation} \label{eq hol}
\pd_s u + J(u) \pd_t u = 0 , \qquad u(s,i)\in L_{H_i} \quad\forall s\in\R,i=0,1 .
\end{equation}
Here we choose the almost complex structure $J$ on $\cR_\S$ that 
is induced by the Hodge operator of some fixed metric $g_\S$ on~$\S$.
Let us first assume that $u$ takes values in the irreducible representations, 
so the pseudoholomorphic equation for $u$ 
actually makes sense since $\cR_\S$ is smooth near its image.
If we consider a lift ${A:\R\times[0,1]\to\cA(\S)}$ of $u$, then this means that
every $A(s,t)$ has stabilizer $\{\pm\one\}\subset\cG(\S)$,
or equivalently $\rd_{A(s,t)}$ is injective on $\Om^0(\S;\su(2))$.
This lift is not unique, but it always takes values in $\Af(\S)$.
So for every $A=A(s,t)$ one has the Hodge decomposition~(\ref{Hodge})
$$
\Om^1(\S;\su(2))=\rd_A\Om^0(\S;\su(2)) \oplus *\rd_A\Om^0(\S;\su(2)) \oplus h^1_A.
$$
Here $h^1_A=\ker\rd_A\cap\ker\rd_A^*\cong \rT_{[A]}\cR_\S$
and $\rd_A\Om^0(\S;\su(2))$ is the tangent space of the $\cG(\S)$-orbit through $A$.
So one can express ${\pd_s u + J(u) \pd_t u = 0}$
in terms of the lift:
The projection of ${\pd_s A + *\pd_t A}$ onto $h^1_A \cong \rT_{[A]}\cR_\S$ vanishes;
i.e.\ $\pd_s A  + *\pd_t A = \rd_A\P + * \rd_A\Psi$
for some ${\P,\Psi:\R\times[0,1]\to\Om^0(\Si;\su(2))}$.
More precisely, (\ref{eq hol}) for $u$ mapping to the irreducible representations
is equivalent to the existence of a lift ${A:\R\times[0,1]\to\cA(\S)}$, 
$u(s,t)=[A(s,t)]$, 
and some ${\P,\Psi:\R\times[0,1]\to\Om^0(\S;\su(2))}$ such that
\begin{equation}\left\{   \label{eq0}
\begin{aligned}
\pd_s A -\rd_A\P + *\bigl( \pd_t A - \rd_A\Psi \bigr) &= 0, \\
*F_A &=0, \\
A(s,i)\in \cL_{H_i} \quad\forall s\in\R,i&= 0,1 .
\end{aligned}\right.
\end{equation}
One can also consider this as a boundary value problem for the connection 
$\P\ds+\Psi\dt+A$ on $\R\times[0,1]\times\S$.
Just note that $A$ determines $\P$ and $\Psi$ uniquely since 
${\laplace_A\P=\rd_A\pd_sA}$, ${\laplace_A\Psi=\rd_A\pd_tA}$,
and $\laplace_A=\rd_A^*\rd_A$ is invertible for irreducible $A=A(s,t)$.
If $A$ is allowed to become reducible, then $\P$ and $\Psi$ have some
extra freedom. If for example $A\equiv 0$, then any two functions
$\P,\Psi:\R\times[0,1]\to\su(2)$ would provide a solution of (\ref{eq0}).
Quotienting out by the gauge action, this moduli space is still infinite dimensional.
We expect however that one can use perturbations of (\ref{eq0}) 
to obtain finite dimensional smooth moduli spaces of trajectories in the cases
that are relevant for $\HF^{\rm symp}_*(\cR_\S,L_{H_0},L_{H_1})$, 
i.e.\ when at least one critical point is irreducible.
Once this symplectic Floer homology is defined via (\ref{eq0}), one should be able
to adapt the adiabatic limit in \cite{DS} to this boundary value problem and 
establish the following second step towards the Atiyah-Floer conjecture.

\begin{con}\label{con:02}
If $Y=H_0\cup_\Sigma H_1$ is a Heegard splitting of a homology 3-sphere,
then there is a natural isomorphism
$$
    \HF^{\rm inst}_*([0,1]\times\Si,\cL_{H_0}\times \cL_{H_1})
    \cong \HF^{\rm symp}_*(\cR_\Si,\cL_{H_0},\cL_{H_1}).
$$
\end{con}

Again, the critical points of both Floer theories are naturally identified.
In the instanton Floer homology the critical points are flat connections on 
$A+\Psi\dt$ on $[0,1]\times\Si$ 
(where flatness means $F_A=0$ and $\dot A - \rd_A\Psi = 0$)
with boundary values $A(0)\in\cL_{H_0}$ and $A(1)\in\cL_{H_1}$.
One can always make $\Psi$ vanish by a gauge transformation, then 
$A$ becomes $t$-independent, so $A(0)=A(1)\in\cL_{H_0}\cap\cL_{H_1}$.
Thus the gauge equivalence classes of these critical points can be identified with 
intersection points of the Lagrangian submanifolds $L_{H_0}$ and $L_{H_1}$ in the 
moduli space $\cR_\Si$ -- which are exactly the
critical points of the symplectic Floer homology.

In order to identify the moduli spaces of trajectories we can choose an appropriate
metric on $[0,1]\times\S$ in the definition of the instanton Floer homology. 
Let us fix the metric $g_\S$ on $\S$ as in (\ref{eq0}) and consider the family 
of metrics $\dt^2 + \ep^2 g_\S$ for $\ep>0$.
With respect to these metrics the trajectory equation (\ref{eq 1}) 
of the instanton Floer homology becomes
\begin{equation}
\label{eqep}
\left\{
\begin{aligned}
\pd_s A -\rd_A\P + *\bigl(\pd_t A - * \rd_A\Psi\bigr) &= 0, \\
\pd_s\Psi - \pd_t\P + [\P,\Psi] + \ep^{-2} * F_A & = 0 , \\
 A(s,i) \in \cL_{H_i} \qquad \forall s\in\R, i&\in\{0,1\} ,
\end{aligned}\right.
\end{equation}
for the triple of ${A:\R\times[0,1]\to\cA(\S)}$ 
and ${\P,\Psi:\R\times[0,1]\to\Om^0(\S;\su(2))}$.
Their energy 
$$
\cE(A,\Phi,\Psi) =
\int_{\R\times[0,1]\times\Sigma} |\pd_s A - \rd_A\Phi|^2 
+ \ep^{-2} |F_A|^2
$$
is determined, independently of $\ep$, by the index 
and the limits at $\pm\infty$ (via a monotonicity formula).
Analogously to \cite{DS} we expect that sequences of such anti-self-dual instantons 
for $\ep\to 0$ converge (modulo gauge) to solutions of (\ref{eq0}). 
Now the gauge equivalence classes of these solutions would exactly be the 
trajectories of the symplectic Floer homology.
Conversely, an implicit function argument should show that for sufficiently small $\ep>0$
near every solution of (\ref{eq0}) one finds a solution of (\ref{eqep}).
This would give the required bijection between the trajectories of the symplectic 
and the instanton Floer homology.

Dostoglou and Salamon indeed dealt with the same equations.
However, they considered a mapping torus $\R\times\S/\sim$ 
(with $(t+1,z)\sim(t,f(z))$ for some diffeomorphism $f$ of $\S$)
instead of our manifold with boundary ${[0,1]\times\S}$, so the boundary conditions in 
(\ref{eq0}) and (\ref{eqep}) are replaced by a twisting condition.
The analytic setup for the definition of the new instanton Floer homology should
also allow to deal with the boundary conditions in this context.
There are however additional difficulties due to reducible connections on the trivial
$\SU(2)$-bundle over $\S$, whereas \cite{DS} deals with the nontrivial 
$\SO(3)$-bundle over $\S$ that has no reducible connections.

\section{Lagrangians in the space of connections}
\label{sec:L}

The purpose of this section is to describe some more properties of
the Lagrangian submanifolds in the space of connections
that were introduced in section~\ref{sec:gauge}.
We again consider more generally a trivial $\rG$-bundle over a Riemann surface $\S$, 
where $\rG$ is any compact Lie group with Lie algebra $\cg$.
We fix $p>2$, then the space of $L^p$-regular connections
$\cA^{0,p}(\S)$ is a symplectic Banach space with symplectic form $\o$ given by
(\ref{omega}).
The gauge group $\cG^{1,p}(\S)$ acts smoothly on $\cA^{0,p}(\S)$ and preserves $\o$.
Moreover, recall that if we equip $\S$ with any Riemannian metric, 
then the corresponding Hodge $*$ operator induces 
an $\o$-compatible complex structure on $\cA^{0,p}(\S)$.

We have proven in~\cite[Theorem~3.1]{W Cauchy} that an $L^p$-connection is flat
in the weak sense iff it is gauge equivalent to a smooth flat connection.
So for our purposes here we simply define the space of flat $L^p$-connections as
$\cA^{0,p}_{\rm flat}(\S):=\cG^{1,p}(\S)^*\Af(\S)\subset\cA^{0,p}(\S)$.
With this definition it is clear that the 
{\bf based holonomy} at any $z\in\S$ is well-defined as a map
$$
\hol_z:\cA^{0,p}_{\rm flat}(\Si)\to{\rm Hom}(\pi_1(\Si),\rG) .
$$
(Here and in the following one actually has to fix one point $z$ in each 
connected component of $\S$.)
It is invariant under the {\bf based gauge group}
$$
\cG^{1,p}_z(\S) := \bigl\{ u\in\cG^{1,p}(\S) \st u(z)=\one \bigr\} .
$$
Next, we call a Banach submanifold $\cL\subset\cA^{0,p}(\S)$ {\bf Lagrangian} 
if is isotropic, $\o|_\cL\equiv 0$,
and coisotropic in the sense of the following implication for all $A\in\cL$ and
$\a\in\cA^{0,p}(\S)$:
If $\o(\a,\b)=0$ for all $\b\in\rT_A\cL$, then $\a\in\rT_A\cL$.
The main properties of gauge invariant Lagrangian submanifolds are summarized 
below. For proofs see \cite[Lemma~4.2,4.3]{W Cauchy}.
(In the case $\rG=\SU(2)$ and for any other connected, simply conected Lie
group with discrete center, the gauge invariance and Lagrangian property
imply that $\cL$ lies in the flat connections; for general groups we make
this additional assumption.)

\begin{lem}  \label{Lagrangian lemma}
Let $\cL\subset\cA^{0,p}(\Si)$ be a Lagrangian submanifold.
Suppose that $\cL\subset\cA^{0,p}_{\rm flat}(\Si)$ and that $\cL$ is invariant 
under the action of $\cG^{1,p}(\Si)$. Then the following holds:
\begin{enumerate}

\item
$\cL\subset(\cA(\S),*)$ is totally real with respect to the Hodge $*$ operator for any
metric on~$\S$. That is $\Om^1(\S;\su(2))=\rT_A\cL\oplus * \rT_A\cL$ for all $A\in\cL$.

\item 
Fix any $z\in\Si$. 
Then $\cL$ has the structure of a principal $\cG^{1,p}_z(\Si)$-bundle
$$
\cG^{1,p}_z(\Si) 
\hookrightarrow  \cL \overset{\hol_z}{\longrightarrow} M .
$$
Here $M\subset{\rm Hom}(\pi_1(\Si),\rG)$ is a smooth manifold of dimension 
$g\cdot\dim \rG$.

\end{enumerate}
\end{lem}

Property~(i) is crucial for the elliptic theory for the boundary value problem 
(\ref{new flow}) in the proof of theorem~\ref{thm:C}. Property~(ii) gives rise to
Banach submanifold coordinates for the Lagrangian that fit well with the Hodge
decomposition of $\Om^1(\S;\su(2))$. This also is the crucial point that forces us
to work on $L^p$-spaces with $p>2$. One does not have a corresponding statement
for Lagrangians in $\cA^{0,2}(\S)$ unless one can find a generalization of the
based gauge group in the $W^{1,2}$-regular gauge transformations.
This would have to be a subgroup that acts freely but has finite codimension.

Next, we consider the Lagrangians given by handle bodies. 
For that purpose we suppose that $\rG$ is connected and simply connected 
and that $\S=\pd H$ is the boundary of a handle body $H$. 
(Both $H$ and $\S$ might have several connected components, in which case
'fixing $z\in\S$' below should be replaced by 'fixing a point in each component'.)

Let $\cL_H$ be the $L^p(\Si)$-closure of the set of smooth flat 
connections on $\Si$ that can be extended to a flat connection on $H$,
$$
\cL_H \,:=\; {\rm cl}\, \bigl\{ A\in\cA_{\rm flat}(\Si) \st \exists 
                              \tA\in\cA_{\rm flat}(H) : \tA|_\Si=A \bigr\} 
\;\subset\;\cA^{0,p}(\Si).
$$
Here again the assumption $p>2$ is crucial for the subsequent properties.
In particular, it is not clear whether the $L^2$-closure is a smooth submanifold.

\begin{lem} \hspace{1mm} \label{LY} \cite[Lemma 4.6]{W Cauchy}\\
\vspace{-5mm}
\begin{enumerate}
\item
$\cL_H = 
 \bigl\{ u^*(A|_\Si) \st A\in\cA_{\rm flat}(H), u\in\cG^{1,p}(\Si) \bigr\} $
\item  $\cL_H\subset\cA^{0,p}(\Si)$ is a Lagrangian submanifold.
\item  $\cL_H\subset\cA^{0,p}_{\rm flat}(\Si)$ and $\cL_H$ is invariant under 
       the action of $\cG^{1,p}(\Si)$.
\item  Fix any $z\in\Si$. Then
\footnote{Here we identify
$
{\rm Hom}(\pi_1(H),\rG) \cong 
\bigl\{ \rho\in{\rm Hom}(\pi_1(\Si),\rG) \st 
\rho(\pd \pi_2(H,\Sigma))=\{\one\}\bigr\}
$.
}
$$
\cL_H = \bigl\{ A\in\cA^{0,p}_{\rm flat}(\Si) \st 
    {\rm hol}_z(A) \in {\rm Hom}(\pi_1(H),\rG) \subset {\rm Hom}(\pi_1(\Si),\rG)  \bigr\} ,
$$
So $\cL_H$ obtains the structure of a $\cG^{1,p}_z(\Si)$-bundle 
over the $g$-fold product 
$M=\rG\times \cdots \times \rG \cong {\rm Hom}(\pi_1(H),\rG)$,
$$
\cG^{1,p}_z(\Si) 
\hookrightarrow  \cL_H
\overset{{\rm hol}_z}{\longrightarrow} {\rm Hom}(\pi_1(H),\rG) .
$$
\end{enumerate}
\end{lem}

%
%
%

Next, although the Lagrangian $\cL_H$ does not necessarily have a smooth
$L^2$-closure, the $L^2(\S)$-norm on $\cL_H$ can be used to control the
corresponding flat connections on $H$ in $L^3(H)$.
This extension property is the crucial trick that circumvents dealing with
the $W^{1,2}$-topology on the gauge group.

\begin{lem} \hspace{1mm}  \label{lem:ext}
There exists a constant $C_H$ such that the following holds.
\begin{enumerate}
\item 
For every smooth path $A : (-\ep,\ep)\to \cL_H\cap\cA(\S)$ there exists a path
$\tA : (-\ep,\ep)\to \cA_{\rm flat}(H)$ with $\tA(s)|_{\pd H}=A(s)$ such that 
$$
\|\pd_s\tA(0)\|_{L^3(H)} \leq C_H \|\pd_s A(0)\|_{L^2(\S)} .
$$
\item
For all $A_0, A_1\in\cL_H\cap\cA(\S)$ there exist extensions
$\tA_0,\tA_1\in\cA_{\rm flat}(H)$ with $A_i=\tA_i|_{\pd H}$ 
such that
\begin{equation} \label{b est}
\|\tA_0-\tA_1\|_{L^3(H)} \leq C_H \|A_0-A_1\|_{L^2(\S)} .
\end{equation}
\end{enumerate}
\end{lem}

The proof in \cite[Lemma 1.6]{W bubbling} uses the coordinates in lemma~\ref{LY}~(iv).
Extensions with the correct holonomy can be constructed by hand, and the estimates
are immediate on this finite dimensional part. For dealing with the
gauge transformations the crucial fact is that there is a continuous extension
operator from $W^{1,2}(\S)$ to $W^{1,3}(H)$.
In (i) this fact is used for functions with values in $\su(2)$, whereas (ii) requires
the nonlinear version for maps to $\SU(2)$. The latter is a nontrivial construction
of Hardt-Lin \cite{HrL} in this borderline Sobolev case (the maps are not automatically
continuous).

\section{Rough guide to the analysis}
\label{sec:analysis}

In this section we give outlines of the proofs of theorems~\ref{thm:C},~\ref{thm:EQ}, and
~\ref{thm:RemSing} for instantons with Lagrangian boundary conditions.\footnote{
The methods will be suitable for generalization to gauge invariant Lagrangians 
as on page~\pageref{Lag ass}. 
The special form of the Lagrangians arising from handle bodies
is only used for the bound on $\frac\pd{\pd\n}e$ in lemma~\ref{lem:Lap} 
and for the isoperimetric inequality in proposition~\ref{thm B}.
}
The detailed proofs can be found in \cite{W elliptic, W bubbling}.
They actually hold for more general domains and metrics than considered here, which 
becomes important when proving the metric independence of the Floer homology,
and when defining products.
We study the boundary value problem (\ref{new flow})
for $\SU(2)$-connections $\X\in\cA(\H^2\times\Si)$,
\begin{equation}\label{bvp}
F_\X + *F_\X = 0,\qquad
\X|_{\{(s,0)\}\times\Si} \in\cL_H \quad\forall s\in\R .
\end{equation}
Here $\H^2=\{(s,t)\in\R^2\st t\geq 0\}$ denotes the half space and we equip $\H^2\times\S$
with a metric $\ds^2 + \dt^2 + g_{s,t}$,
where the metric $g_{s,t}$ on $\S$ varies smoothly with
$(s,t)\in \H^2$ and is constant outside of a compact subset.

\subsection{Proof of Compactness} \label{sec:C} \hspace{2mm}

For all results in this subsection 
the Lagrangian $\cL_H$ in (\ref{bvp}) can be replaced by a general
gauge invariant Lagrangian submanifold $\cL\subset\cA^{0,p}(\S)$.
The compactness theorem~\ref{thm:C} in case (in-L) is a consequence
of the following lemma and theorem.
The lemma yields the local $L^p$-bounds that are assumed in the theorem.
It is based on mean value inequalities and will thus be proven later in section~\ref{mean}.
Here $B_r(x)\subset\R^2$ is the closed $2$-dimensional ball 
of radius $r>0$ centered at $x\in\R^2$, and we denote
$D_r(x):=B_r(x)\cap\H^2$.
In particular, $D_r:=D_r(0)\subset\H^2$ is the closed half ball of radius $r$.

\begin{lem} \label{lem L2 bound} \cite[Lemma 2.4]{W bubbling}  
Let $\X^\n\in\cA(\H^2\times\S)$ be a sequence of anti-self-dual connections
and suppose that for some $x_0\in\H^2$ and $\d>0$
$$
\sup_{\n} \sup_{x\in D_{2\d}(x_0)} 
\bigl\| F_{\X^\n} (x) \bigr\|_{L^2(\Si)} <\infty .
$$
Then for every $2<p<3$
$\displaystyle\qquad
\sup_{\n} \bigl\| F_{\X^\n} \bigr\|_{L^p(D_\d(x_0)\times\Si)} < \infty .
$\\
If in fact $B_\d(x_0)\cap\pd\H=\emptyset$, then moreover
$\displaystyle\quad
\sup_{\n} \bigl\| F_{\X^\n} \bigr\|_{L^\infty(B_\d(x_0)\times\Si)} < \infty 
$.
\end{lem}

\begin{thm}  \label{thm:cpt}
Let $p>2$.
Suppose that $\X^\n\in \cA(\H^2\times\S)$ is a sequence of solutions of (\ref{bvp})
such that $\sup_\n\|F_{\X^\n}\|_{L^p(K)}<\infty$ for every compact
subset $K\subset\H\times\S$.
Then there exists a subsequence (again denoted by $\X^\n$) and a sequence of gauge transformations 
$u^\n\in\cG(\H^2\times\S)$ such that $u^{\n\;*}\X^\n$ converges uniformly with all derivatives on every
compact subset of $\H^2\times\S$.
\end{thm}

Note that it is crucial to establish this compactness for $2<p<3$ since
the previous lemma only provides those curvature bounds near the boundary.
Next, we outline the steps of the proof in \cite[Theorem B]{W elliptic}
of theorem~\ref{thm:cpt}.
By standard gauge theoretic arguments it boils down to the boundary regularity
theory in 5b)--f) below. The crucial step is f), where the Lagrangian 
enters as totally real boundary condition for a Cauchy-Riemann equation. 
The case $2<p\leq 4$ requires a separate treatment described in a') and~f').

\noindent
{\bf 1) Reduction to compact domains: }
By a Donaldson-Kronheimer trick 
\cite[Prop.~7.6]{W}
it suffices to prove the assertion on $D_k\times\S$ for every $k\in\N$.
Then the gauge transformations on $D_k\times\S$ 
can be extended to $\H^2\times\S$ and can be
interpolated with gauge transformations obtained
on larger domains. A diagonal subsequence
then satisfies the claimed $\cC^\infty_{\rm loc}$-convergence 
on $\H^2\times\S$.

So we consider a sequence $\X^\n\in\cA(\H^2\times\S)$ of solutions 
whose curvature is in particular $L^p$-bounded on $\cU\times\S$, where
$\cU\subset\H^2$ is some compact domain with smooth boundary and $D_k\subset{\rm int}(\cU)$.
Then we need to find gauge transformations and a convergent subsequence on $D_k\times\S$.

In the subsequent steps one frequently gets a new estimate only on a smaller domain
$\cU_i\subset{\rm int}(\cU)$. (Note that the interior includes points on $\pd\H^2$.)
However, we can always choose these such that $D_k\subset{\rm int}(\cU_i)$.

\noindent
{\bf 2) Weak convergence:}
We can apply Uhlenbeck's weak compactness theorem~\ref{thm:weak comp} on $\cU\times\S$.
It provides a subsequence (still denoted $\X^\n$)
and gauge transformations $u^\n\in\cG^{2,p}(\cU\times\S)$ such that 
${u^{\n\;*}\X^\n \to \X^\infty}$ in the weak $W^{1,p}$-topology
with a limit connection $\X^\infty\in\cA^{1,p}(\cU\times\S)$.

\noindent
{\bf 3) Regularity for limit solution:}
The limit $\X^\infty$ now also solves (\ref{bvp}).
For the boundary conditions this is due to the compact Sobolev embedding 
$W^{1,p}(\cU\times\S)\hookrightarrow\cC^0(\cU,L^p(\S))$. 
In the nonstandard case $2<p\leq 4$ this embedding is established 
in \cite[Lemma~2.5]{W elliptic}.

Now one finds a gauge transformation $u\in\cG(\cU_1\times\S)$ such that
$u^*\X^\infty$ is smooth on $\cU_1\times\S$. (This is proven analogously to
the iteration in 5), with estimates replaced by regularity statements. 
For the local slice theorem in 4) it suffices to pick a smooth connection
$\X_0$ that is $W^{1,p}$-close to $\X=\X^\infty$.)
One thus finds that $(u^\n u)^*\X^\n \to \X_0$ in the weak $W^{1,p}$-topology
on $\cU_1\times\S$, with a smooth limit $\X_0=u^*\X^\infty$.

\noindent
{\bf 4) Relative Coulomb gauge:}
Next, the local slice theorem~\ref{thm:local slice} provides a sequence of gauge transformations
$v^\n\in\cG(\cU_1\times\S)$ such that still $v^{\n\;*}\X^\n \to \X_0$ converges $W^{1,p}$-weakly,
and in addition each $\X=v^{\n\;*}\X^\n$ satisfies
\begin{equation}\label{relCoul}
\rd_{\X_0}^*(\X-\X_0)=0 , \qquad
*(\X-\X_0)|_{\pd\cU_1\times\S}=0 .
\end{equation}

\noindent
{\bf 5) Elliptic estimates for (\ref{bvp})\&(\ref{relCoul}):}
From 2--4) we have a subsequence and gauge transformations $v^\n$ such that each 
$\X=v^{\n\;*}\X^\n$ satisfies (\ref{bvp}), (\ref{relCoul}), and 
$\|\X\|_{W^{1,p}(\cU_1\times\S)}\leq C_1$ for some uniform constant $C_1$.
By iterating the following steps a)--f) one next finds uniform constants $C_\ell$
such that
$\|\X\|_{W^{\ell,p}(\cU_\ell\times\S)}\leq C_\ell$ for all $\ell\in\N$ and for all
$\X=v^{\n\;*}\X^\n$.
Finally, due to the compact Sobolev embeddings 
$W^{\ell,p}(D_k\times\S)\hookrightarrow\cC^{\ell-2}(D_k\times\S)$ 
one then finds a diagonal subsequence that 
converges with all derivatives on $D_k\times\S$.
This is what was to be shown according to 1).

For a)--f) we give the arguments in the case $\ell=2$ and $p>4$.
This first step is considerably harder for $2<p\leq 4$ and requires 
a separate iteration, which is roughly indicated in a') and f').
The iteration for $\ell\geq 3$ and any $p>2$ then works completely
analogous to the arguments below.

\noindent
{\bf a) Interior estimates:}
From (\ref{bvp}) and (\ref{relCoul}) we obtain the Hodge Laplacian
\begin{align} \label{Lap}
\laplace \X 
&= - \rd^* \bigl( \half[\X\wedge\X] + \half *[\X\wedge\X]  \bigr)   
 + \rd\rd_{\X_0}^*\X_0  + \rd * [\X_0\wedge *\X]  . 
\end{align}
Here the right hand side is bounded in $L^p$, and the leading order of the
left hand side in local coordinates is the Laplacian on the components of $\X$.
Thus the elliptic estimate for the Laplace equation yields a $W^{2,p}$-bound
on $\X$ in the interior of $\cU_1\setminus\pd\H^2$.

Going through the arguments up to this point also proves theorem~\ref{thm:C}
in the instanton case (inst) without boundary.

\noindent
{\bf a') Special iteration for $\mathbf{W^{2,p}}$-bounds with 
$\mathbf{2<p\leq 4}$:}
In this case the right hand side of (\ref{Lap}) lies in $L^q$ for some $q<p$,
so one only obtains a $W^{2,q}$-bound. However, by a Sobolev embedding,
this also gives a $W^{1,p'}$-bound for some $p'>p$. 
Iteration of a) then yields $W^{2,q_i}$-bounds for a strictly increasing
sequence which reaches $q_N\geq p$ after finitely many steps.
An analogous iteration will work for steps b--f).

\noindent
{\bf b) Splitting the equation near the boundary:}
It remains to obtain a $W^{2,p}$-bound on $\X$ near $D_k\cap\pd\H^2$.
For that purpose we rewrite (\ref{bvp}) and (\ref{relCoul})
in the splitting $\X=\P\ds +\Psi\dt + A$ (and analogous for the smooth $\X_0$) 
with $\P,\Psi\in W^{1,p}(\cU_1\times\S;\su(2))$
and $A\in W^{1,p}(\cU_1\times\S,\rT^*\S\otimes\su(2))$,
\begin{equation} \label{bbvp} 
\left\{ \begin{aligned}
( \pd_s A - \rd_A\P ) + * ( \pd_t A - \rd_A\Psi ) & = 0 ,\\
\pd_s \Psi - \pd_t\P + [\P,\Psi] + * F_A  & = 0 , \\
\nabla_s(\P-\P_0) + \nabla_t(\Psi-\Psi_0) - \rd_{A_0}^*(A-A_0) & = 0 , \\
( \Psi - \Psi_0 )|_{t=0} & = 0 , \\
 A(s,0) \in \cL \quad &\forall s\in\R .
\end{aligned} \right.
\end{equation}
Here we use the notation $\nabla_s = \pd_s + [\P_0,\cdot]$ and 
$\nabla_t= \pd_t + [\Psi_0,\cdot]$.

\noindent
{\bf c) Estimates for $\mathbf\Psi$:}
From (\ref{Lap}) we know that $\laplace\Psi$ is $L^p$-bounded.
In addition, we have the inhomogeneous Dirichlet condition 
$ \Psi|_{t=0} = \Psi_0 |_{t=0} $.
Thus the elliptic estimate for the Dirichlet boundary value problem
implies a $W^{2,p}$-bound on $\Psi$ up to the boundary.

\noindent
{\bf d) Estimates for $\mathbf\P$:}
Again, $\laplace\P$ is $L^p$-bounded due to (\ref{Lap}).
Moreover, we have an inhomogeneous Neumann condition 
$\pd_t\P|_{t=0} = (\pd_s \Psi_0  + [\P,\Psi_0])|_{t=0}$
since the Lagrangian boundary condition with $\cL\subset\cA_{\rm flat}^{0,p}(\S)$
in particular implies $F_A|_{t=0}=0$.
Then the elliptic estimate for the Neumann boundary value problem
(e.g.\ \cite[Theorem~3.2]{W}) provides a $W^{2,p}$-bound on $\P$.

\noindent
{\bf e) Estimates for $\mathbf{\nabla_\S A}$:}
We can now rewrite (\ref{bbvp}) to
express the differential and codifferential of $A(s,t)\in\Om^1(\S;\su(2))$
for every $(s,t)\in\cU_1$ as 
\begin{align*}
* \rd_\S A &= -\half * [A\wedge A] - \pd_s \Psi + \pd_t\P - [\P,\Psi] , \\
\rd_\S^* A &= \nabla_s(\P-\P_0) + \nabla_t(\Psi-\Psi_0) + * [A_0\wedge*A] - \rd_{A_0}^*A_0 .
\end{align*}
Due to the previously established bounds on $\P$ and $\Psi$ the right hand sides 
are bounded in $W^{1,p}(\cU_2\times\S)$, that is in
$W^{1,p}(\cU_2,L^p(\S))$ and $L^p(\cU_2,W^{1,p}(\S))$.
Now the elliptic estimates from the Hodge decomposition for each $(s,t)\in\cU_2$
can be integrated to give bounds on $\nabla_\S A$ in the same spaces, and hence
a $W^{1,p}$-bound. For a detailed statement and proof see \cite[Lemma~2.9]{W elliptic}.

\noindent
{\bf f) Estimates for $\mathbf{\pd_s A, \pd_t A}$:}
So far $A$ is bounded in $L^p(\cU_2,W^{2,p}(\S))$ and $W^{1,p}(\cU_2,W^{1,p}(\S))$.
To achieve a $W^{2,p}$-bound it remains to find an estimate in $W^{2,p}(\cU_2,L^p(\S))$,
that is on $\pd_s A$ and $\pd_t A$.
At this point, the full Lagrangian boundary condition needs to be used.
Up to now, we only used its local part, the slice-wise flatness.
The additional holonomy conditions are of global type 
(requiring knowledge of the connection on loops in $\S$),
so this information is lost when one localizes, i.e.\ goes to a coordinate chart in $\S$.

The solution is to consider $A$ as map from $\cU_1$ to the Banach space $\cA^{0,p}(\S)$.
This is a complex space when equipped with the Hodge $*$ operator.
So we can rewrite (\ref{bbvp}) and recall lemma~\ref{Lagrangian lemma}~(i) 
to see that $A$ satisfies a Cauchy-Riemann equation
with totally real boundary conditions:
\begin{equation} \label{holo}
\pd_s A  + * \pd_t A   = \rd_A\P + *\rd_A\Psi , \qquad
A(s,0) \in \cL \quad \forall s\in\R .
\end{equation}
Now one basically has to go through the proof of theorem~\ref{thm:C} for the
holomorphic curves in case (symp) with the extra difficulty that the target space
is infinite dimensional. This would be fairly standard for a Hilbert space.
However, the iteration only works for $p>2$ and we also need to work with $p>2$ 
to make sure that the Lagrangians are smooth submanifolds.

Here we use the general theory in \cite{W Cauchy}
for maps to a complex Banach space $X$. 
The crucial assumption is that $X$ is a closed subspace of an 
$L^p$-space for some ${1<p<\infty}$ on a closed manifold (for example $X=\cA^{0,p}(\S)$).
Then the elliptic $L^p$-estimates (with the same Sobolev exponent as in $X$) 
hold for the Dirichlet and Neumann problem.
One can then use the usual argument for the Cauchy-Riemann 
equation with totally real boundary conditions:
In a submanifold chart the components of the map 
${u:\cU\to \rT_{z_0}\cL\times \rT_{z_0}\cL\cong X}$ 
satisfy Dirichlet and Neumann boundary conditions --
at the expense of the complex structure becoming $u$-dependent. 
From an $L^p$-bound on ${(\pd_s^2 +\pd_t^2)u}$ one then obtains a $W^{2,p}$-estimate on $u$.
Due to the nonlinearity in the complex structure however,
the $L^p$-estimate on 
${(\pd_s^2 +\pd_t^2)u}$
requires $W^{1,2p}$-bounds on $u$ and $\pd_s u + J\pd_t u$.

In (\ref{holo}) the right hand side is bounded in 
$W^{1,p}(\cU,L^p(\S))$ due to the previous bounds on $\P$ and $\Psi$
(on some domain $\cU \subset {\rm int}(\cU_1)$ with ${\cU_2\subset{\rm int}(\cU)}$).
By the above discussion we now have to write $p=2p'$ and we only obtain 
$W^{2,p'}$-estimates for $A:\cU\to X=\cA^{0,p'}(\S)$ 
with $\pd_s A + *\pd_t A \in W^{1,2p'}(\cU,X)$.
So this last step yields a bound on $\X$ in $W^{2,\frac p2}(\cU_2\times\S)$.
For $p>4$ we still have $\frac p2>2$ and the further iteration
yields $W^{\ell,\frac p2}$-bounds for all $\ell\in\N$.

\noindent
{\bf f') Special case $\mathbf{2<p\leq 4}$ for $\mathbf{W^{2,p}}$-bounds:}
In this case we only have $q<p$ in the $W^{2,q}$- and $W^{1,q}$-bounds 
on $\P,\Psi$, and $\nabla_\S A$ from c)--e). 
So the right hand side in (\ref{holo}) is of even lower regularity that will not
fit in the above arguments. However, it is bounded in $L^r(\cU,L^p(\S))$ for some
$r>p$. So one can use the submanifold charts for $\cL\subset\cA^{0,p}(\S)$
to write $A$ as a map $u:\cU\to \rT_{A_0}\cL\times \rT_{A_0}\cL$, 
where $\rT_{A_0}\cL\subset \cA^{0,p}(\S)$ is a closed subspace.
The two components of $u$ then satisfy weak Dirichlet and Neumann equations 
with the weak Laplacian in $W^{-1,r}(\cU,L^p(\S))$.
The previous general theory unfortunately only works when we replace the $r>p$
by $p$ and it would then give a bound on $u$ in $W^{1,p}(\cU_2,L^p(\S))$, 
which is what we started out with.
However, one can use all the usual elliptic estimates when the target is
a Hilbert space. So we consider $u$ as map into $\cA^{0,2}(\S)\times\cA^{0,2}(\S)$
with a $W^{-1,r}$-bound on its weak Laplacian.
This yields a $W^{1,r}(\cU_2,L^2(\S))$-bound on $u$ with $r>p$.
The previous bounds in e) moreover imply a $W^{1,q}(\cU_2,L^s(\S))$-bound on $u$,
where $q<p$ but $s>p$ since it results from the Sobolev embedding
$W^{1,q}(\S)\hookrightarrow L^s(\S)$.
Now these two bounds can be interpolated to obtain a 
$W^{1,p'}(\cU_2\times\S)$-bound with $p'>p$. This bound on $u$ also translates 
into a $W^{1,p'}$-bound on $A$, which fits into the same iteration as in a').

\medskip

\subsection{Mean value inequalities} \label{mean} \hspace{2mm}

The proof of theorems~\ref{thm:EQ} and \ref{thm:RemSing} as well as
lemma~\ref{lem L2 bound} makes use of some mean value inequalities which
we summarize here.
These are based on a generalization of the mean value inequality 
for subharmonic functions. Here we state it for the Euclidean half space
$\H^n$. In the interior case this is wellknown for general metrics. In the
case of balls intersecting the boundary this was proven in 
\cite[Theorem~1.3]{W mean} for the Euclidean metric.

\begin{prp}
\label{prp:mean}
For every $n\geq 2$ there exists a constant $C$ and for all $a,b\geq 0$ there exists 
$\hbar(a,b)>0$ such that the following holds: 

Let $D_r(y)\subset\H^n$ be the partial ball of radius $r>0$ and centre $y\in\H^n$. 
Suppose that $e\in\cC^2(D_r(y),[0,\infty))$ satisfies for some constants $A,B\geq 0$
\begin{equation*}
\left\{\begin{array}{ll}
\laplace e &\leq  A e + a e^{\frac {n+2}n} , \\
\frac\pd{\pd\n}\bigr|_{\pd\H^n} e \hspace{-2mm} &\leq B e + b e^{\frac {n+1}n} ,
\end{array}\right.
\qquad\text{and}\qquad
\int_{D_r(y)} e \leq \hbar(a,b) .
\end{equation*}
Then
$\displaystyle\quad
e(y) \leq   C \bigl(A^{\frac n2} + B^n + r^{-n} \bigr) \tint_{D_r(y)} e .
$
\end{prp}

%
%

For all three types of Floer theory that are discussed in section~\ref{sec:FH},
the energy densities satisfy the differential inequalities for 
proposition~\ref{prp:mean} with exactly the critical nonlinearities.
These estimates are summarized below.

\begin{lem} \label{lem:Lap}
Consider a solution of the trajectory equation (T) in 
definition~\ref{def:trajectory}.
Its energy density $e$ satisfies the following nonlinear
bounds on $\laplace e$ and $\frac\pd{\pd\n} e$
with constants $a,b,C$.

\smallskip

\noindent
{\bf (inst):}
$e=|\pd_s B|^2 :\R\times Y \to [0,\infty)$ satisfies
$\displaystyle\qquad
\laplace e \leq C\, e + a \, e^{\frac 32} $.

\smallskip

\noindent
{\bf (in-L,interior):}
$e =|\pd_s A|^2 + |F_A|^2 :\R\times[0,1]\times\S \to [0,\infty)$ satisfies
\begin{equation*}
\laplace e \leq C\, e + a \, e^{\frac 32} .
\end{equation*}

\smallskip

\noindent
{\bf (in-L,boundary):}
$e=\|\pd_s A\|_{L^2(\Sigma)}^2 +\|F_{A}\|_{L^2(\Sigma)}^2 : 
\R\times[0,1]\to[0,\infty)$ satisfies
\begin{equation*}
\laplace e \leq  C \bigl( 1 + \|F_A\|_{L^\infty(\S)} \bigr) e , \qquad
\tfrac\pd{\pd\n} e \leq C\, e +  b \, e^{\frac 32} .
\end{equation*}

\smallskip

\noindent
{\bf (symp):}
$e=|\pd_s u|^2 : \R\times[0,1]\to[0,\infty)$ satisfies
\begin{equation*}
\laplace e \leq  a \, e^2 , \qquad
\tfrac\pd{\pd\n} e \leq b \, e^{\frac 32} .
\end{equation*}

\end{lem}

\noindent
{\bf Indications of proofs of lemma~\ref{lem:Lap}:}

For the holomorphic curves in case (symp) one picks up linear terms
in the estimates if the almost complex structure $J$ varies over the domain. 
The bound on the Laplacian can be
found in e.g.\ \cite[Lemma 4.3.1]{MS2}. The bound on the normal derivative
was wellknown and is proven in \cite[Lemma A.1]{W mean} using 
Darboux-Weinstein coordinates near the Lagrangian.

For the anti-self-dual instantons in case (inst) this estimate is a direct
consequence of a Bochner-Weitzenb\"ock formula, see e.g.\ 
\cite[Lemma A.2]{W mean}. It was used by Uhlenbeck
\cite[Lemma 3.1]{U1} in a slightly different formulation.
For the anti-self-dual instantons with Lagrangian boundary conditions, one
has the same bound on the Laplacian, as stated in (in-L,interior). 
However, this only provides estimates in the interior
(on balls that do not intersect the boundary) since one does not have a 
bound on the normal derivative.

In view of the global methods in section~\ref{sec:C}~f)
that were necessary for the proof of the basic compactness theorem~\ref{thm:C}
it should not be surprising that we were not able to obtain any bound on
$\frac\pd{\pd\n} e$ in terms of $e$, let alone by $b\,e^{\frac 54}$.
It is highly unclear how the (nonlocal) holonomy part of the 
Lagrangian boundary condition should be utilised for such a local estimate.
On the other hand, there are examples showing that such an estimate cannot
follow only from the (local) flatness part of the Lagrangian boundary condition.

Thus it seems natural that the full Lagrangian boundary condition is only
captured by the 2-dimensional energy density given in (in-L,boundary). 
Indeed, we obtain the same bound on the normal derivative as in case (symp).
The proof in \cite[Lemma~2.3]{W bubbling} works as follows:
A simple calculation using the trajectory equation (T) in 
definition~\ref{def:trajectory} (that is (\ref{bvp}) in temporal gauge)
gives the normal derivative at the $t=0$ boundary component:
\begin{align*}
- \half \tfrac\pd{\pd t} e \bigr|_{t=0}
%
%
&= - \int_\S \la \pd_s A \wedge * \pd_s \bigl( * \pd_s A \bigr) \ra\bigr|_{t=0} \\
&\leq \Bigl( C \bigl\| \pd_s A \bigr\|_{L^2(\S)}^2 
     + \int_\Si \la \pd_s A \wedge \pd_s^2 A \ra \Bigr)\Bigr|_{t=0} .
\end{align*}
%
%
%
%
%
%
%
Recall that $e=\|\pd_s A\|_{L^2(\Sigma)}^2 +\|F_{A}\|_{L^2(\Sigma)}^2$ 
and $F_A\bigr|_{t=0}=0$ by the boundary condition.
So the first term on the right hand side is just $C e$ for a constant~$C$. 
The crucial second term is $\o(\pd_s A,\pd_s^2 A)$ 
for a path ${A:(-\ep,\ep)\to\cL_H}$ in the Lagrangian
and with the symplectic form (\ref{omega}).

This term would vanish if the Lagrangian was straight --
as in Darboux-Weinstein coordinates.
Otherwise the curvature of the Lagrangian leads to a cubic term.
For general infinite dimensional Lagrangians the curvature might not be 
suitably bounded, and it is not clear whether Darboux-Weinstein coordinates
even exist. Fortunately, we are dealing with Lagrangians that are compact 
modulo gauge transformations.
A proof along this line would require a subtle linear estimate
for gauge transformations in the critical Sobolev space $W^{1,2}(\S)$,
which has not been carried out yet.
For the special Lagrangian $\cL_H$ arising from a handle body
we can use the following trick based on the extension property 
in lemma~\ref{lem:ext}~(i).

We have $A(s)=\tA(s)|_{\pd H}$ for a path of extensions 
${\tA : (-\ep,\ep)\to\cA_{\rm flat}(H)}$ such that
${\|\pd_s\tA\|_{L^3(H)} \leq C \|\pd_s A\|_{L^2(\S)}}$.
Now Stokes' theorem gives
\begin{align*}
\int_\S \la \pd_s \tA \wedge \pd_s^2\tA \ra 
&=  \int_H \la \rd_\tA \pd_s \tA \wedge \pd_s^2 \tA \ra  
  - \int_H \la \pd_s \tA \wedge \rd_\tA\pd_s^2 \tA \ra  \\
&=  \int_H \la \pd_s \tA \wedge [\pd_s\tA \wedge \pd_s\tA] \ra  \\
&\leq \|\pd_s\tA\|_{L^3(H)}^3 
\,\leq\, C^3 \|\pd_s A\|_{L^2(\S)}^3
\,=\, C^3 e^{\frac 32} .
\end{align*}
Here we used the fact that $F_\tA\equiv 0$, 
hence $\rd_\tA\pd_s\tA = \pd_s F_\tA = 0$,
and moreover
$0 = \pd_s^2 F_\tA = \rd_\tA\pd_s^2\tA + [\pd_s\tA \wedge \pd_s\tA]$.
This proves 
$\pd e/\pd\n \leq C\, e +  b \, e^{3/2} $.

The price for going to the more global energy density in (in-L,boundary)
has to be paid when considering the Laplacian.
The straight forward calculations in \cite[Lemma~2.3]{W bubbling} yield
$$
\laplace e
\leq C \bigl(\|\pd_s A\|_{L^2(\Sigma)}^2 +\|F_{A}\|_{L^2(\Sigma)}^2 \bigr)
  - 20 \la F_A \,,\, [\pd_s A \wedge \pd_s A] \ra_{L^2(\Si)} .
$$
The first term is just $C e$. The second term should also be bounded in
terms of the $L^2$-norms of the curvature components $\pd_s A$ and $F_A$.
However, the best bound that we can find is 
$\|F_A\|_{L^\infty(\S)} \|\pd_s A\|_{L^2(\Sigma)}^2 
\leq \|F_A\|_{L^\infty(\S)} e $.
Here we use the $L^\infty$-norm on $F_A$ since this has better analytic
properties, in particular Dirichlet boundary conditions $F_A|_{t=0}=0$,
whereas $\pd_s A$ only satisfies Lagrangian boundary conditions
(of global type). 
This will be crucial in the proof of the energy quantization theorem~\ref{thm:EQ},
where we will find that $\laplace e\leq C(1+ \|F_A\|_{L^\infty(\S)}) e$ 
is essentially bounded by $C e^2$.\\

\noindent
{\bf Proof of lemma~\ref{lem L2 bound}:}
This is a consequence of the mean value inequality in 
proposition~\ref{prp:mean} applied to the energy densities 
$e_\n = |\pd_s A^\n|^2 + |F_A^\n|^2 = \half |F_{\X^\n}|^2$
from case (in-L,interior) of lemma~\ref{lem:Lap}.
The assumption can be read as
\begin{equation}\label{bound s}
\int_\S e_\n(x,\cdot) \leq K  
\qquad\text{for all}\; x\in D_{2\d}(x_0) 
\end{equation}
with some uniform constant $K$.
On $4$-dimensional balls $B^4_\ep(y)$ 
that are entirely contained in $D_{2\d}(x_0)\times\S$
this implies $\int_{B^4_\ep(y)} e_\n \leq \pi K \ep^2$.
Now there is a maximal radius $\ep_0\in(0,\d)$ such that
for all $\ep\leq\ep_0$ this energy is less than $\hbar(a)$
and thus one has the mean value inequality
\begin{equation}\label{bound}
e_\n(y)\leq C(1+\ep^{-4})\pi K \ep^2 .
\end{equation}
In the interior case, one fixes a radius $0<\ep\leq\ep_0$ less than the distance
${\rm dist}(B_\d(x_0),\pd\H^2)>0$.
Then all balls $B^4_\ep(y)$ for $y\in B_\d(x_0)\times\S$ are contained
in $D_{2\d}(x_0)\times\S$ and (\ref{bound}) is the claimed uniform bound.

In the boundary case $x_0=(s_0,t_0)$ with $t_0\leq\d$
one cannot use a fixed radius for the balls near the boundary.
At $y=(s,t,z)\in D_\d(x_0)\times\S$ the maximal ball that is 
entirely contained in $D_{2\d}(x_0)\times\S$ has radius $\ep=\min(t,\d)$.
So for all $(s,t,z)\in D_\d(x_0)\times\S$ with $0<t\leq\ep_0$ 
the mean value inequality (\ref{bound}) gives
$$
e_\n(s,t,z)\leq C'(t^2+t^{-2}).
$$
Away from the boundary, for $t\geq \ep_0$, 
this also holds with some modified constant $C'$
by (\ref{bound}) with a fixed radius.
Now this bound blows up as $t\to 0$, but it can be interpolated with
(\ref{bound s}) to give an $L^p$-bound on $|F_{\X^\n}|=(e_\n)^{1/2}$
by the following integral which is finite for $2<p<3$.
$$
\int_{D_\d(x_0)\times\S} (e_\n)^{\frac p2}
\;\leq\; \int_{D_\d(x_0)} \bigl( C'(t^2+t^{-2}) \bigr)^{\frac p2 -1} \int_\S e_\n
\;\leq\; C'' \Bigl( 1 + \int_0^{t_0+\d} t^{2-p} \dt \Bigr) .
$$

\subsection{Proof of Energy Quantization} \hspace{2mm} \label{eq}

The proof of theorem~\ref{thm:EQ} for anti-self-dual instantons without boundary and 
for the holomorphic curves is a direct consequence of the mean value inequality in 
proposition~\ref{prp:mean} applied to the energy densities in lemma~\ref{lem:Lap}.
(See \cite[Theorem~2.1]{W mean} for this general energy quantization principle.)
For (inst) and the interior of (in-L) this is the simplest
version of the argument -- on balls with no boundary condition in dimension~$n=4$.
Here we give the argument for the holomorphic curves in (symp), more generally
for a sequence of energy density functions $e_i:\R\times[0,1]\to[0,\infty)$
satisfying
\begin{equation*}
\laplace e_i \leq K\, e_i + a \, e_i^2 , \qquad\quad
\tfrac\pd{\pd\n} e_i \leq B\, e_i + b \, e_i^{\frac 32} .
\end{equation*}
We need to prove that if the energy densities blow up 
at some $x\in\R\times[0,1]$,
\begin{equation} \label{ass}
\sup_i \; \sup_{D_\d(x)} e_i \;=\; \infty \qquad\forall \d>0 ,
\end{equation}
then (for a subsequence) a fixed energy quantum $\hbar>0$ concentrates there,
\begin{equation} \label{claim}
\int_{D_\d(x)} e_i \;>\; \hbar \qquad\forall \d>0 .
\end{equation}
The same needs to be proven in case (in-L) for boundary points $x$.
For these anti-self-dual instantons with Lagrangian boundary conditions
we use the energy density 
$e_i=\|\pd_s A_i\|_{L^2(\S)}^2 + \|F_{A_i}\|_{L^2(\S)}^2$
as in (in-L,boundary) of lemma~\ref{lem:Lap}.
So the constant $K$ above is replaced by the 
unbounded function $C(1+\|F_{A_i}\|_{L^\infty(\S)})$.
Moreover, the assertion (\ref{claim}) in this case implies the concentration
of energy near $\{x\}\times\S\subset\R\times[0,1]\times\S$.

So let us assume (\ref{ass}).
Then we find a subsequence and points $x_i\to x$ such that
$e_i(x_i)=R_i^2$ blows up with a certain rate $R_i\to\infty$.
We will now try to apply proposition~\ref{prp:mean} on the balls 
$D_{\d_i}(x_i)$ of radius ${\d_i:=R_i^{-1/2}>0}$.
For that purpose we need to assume that
$\int_{D_{\d_i}(x_i)} e_i \leq \hbar=\hbar(a,b)$.
If that is the case then we obtain the mean value inequality
$$
R_i^2 \;=\; e(x_i) 
\;\leq\; C \bigl(K + B^2 + \d_i^{-2} \bigr) \int_{B_{\d_i}(x_i)} e_i .
$$
Multiplication by $R_i^{-2}=\d_i^2R_i^{-1}$ then implies
$$
1 \;\leq\; C\hbar \bigl(K R_i^{-2} + B^2 R_i^{-2} + R_i^{-1} \bigr).
$$
First assume that $K$ is constant. Then the right hand side converges to $0$.
Thus the assumption must have failed for all sufficiently large $i\in\N$,
that is $\int_{D_{\d_i}(x_i)} e_i > \hbar$. This implies
the energy concentration~(\ref{claim}).

If $K$ is not a constant, then this argument still works as long as
$K\leq C'R_i^2$. In that case the limit $i\to\infty$ implies $1 \leq CC'\hbar $. 
If one chooses $\hbar\leq (2 C C')^{-1}$, then this
gives a contradiction and thus proves the energy concentration.

So for anti-self-dual instantons with Lagrangian 
boundary conditions in case (in-L) we have to prove that if 
$e_i=R_i^2$ blows up, then the functions $K=C(1+\|F_{A_i}\|_{L^\infty(\S)})$
are bounded by $C'R_i^2$.
This statement is slightly weaker than a direct bound 
$\|F_A\|_{L^\infty(\S)}\leq C \|F_\X\|_{L^2(\S)}^2 = C e$ would be,
but it still shows that $\laplace e\leq C(1+ \|F_A\|_{L^\infty(\S)}) e$ 
is essentially bounded by $C e^2$.

By using the Hofer trick \cite[6.4 Lemma 5]{HZ} within the previous argument
one can additionally control $e_i$ by the blowup rate 
on small neighbourhoods.
One then needs to establish the following 
as in \cite[Proposition~2.7]{W bubbling}. 

\pagebreak

\begin{lem} {\bf (Crucial Estimate):} \label{crucial}
Let $\X_i=\P_i\ds+\Psi_i\dt+A_i\in\cA(\H^2\times\S)$ 
be a sequence of solutions of (\ref{bvp}). 
Consider a sequence of blowup points $\H^2\ni x_i\to 0$ with the blowup speed 
$R_i\to\infty$. Assume an $L^2(\S)$-control on the full curvature
on (partial) balls of radius $2\ep_i\to 0$ such that $\ep_i R_i\to \infty$,
$$
\|F_{\X_i}(x,\cdot)\|_{L^2(\S)} \leq R_i \qquad\forall x\in D_{2\ep_i}(x_i). 
$$
Then one obtains an $L^\infty$-control on the curvature component 
$$
\|F_{A_i}(x,\cdot)\|_{L^\infty(\S)} \leq C R_i^2 \qquad\forall x\in D_{\ep_i}(x_i) .
$$ 
\end{lem}

The proof combines all previous techniques to a subtle contradiction. 
This is what remains of the usual energy quantization proof via local rescaling:

\noindent
{\bf 1.)~Assume the contrary:}
Then one finds sequences of solutions $\X_i$,
points $(x_i,z_i)\to (0,z)\in\H^2\times\S$, and 
$R_i\to\infty$, $\ep_i\to 0$, $C_i\to\infty$
with $\ep_i R_i\to\infty$,
$$ \textstyle
\sup_{x\in D_{\ep_i}(x_i)} 
\|F_{\X_i}(x,\cdot)\|_{L^2(\S)} \leq R_i , \qquad
|F_{A_i}(x_i,z_i)| \geq (C_i R_i)^2 .
$$ 

\noindent
{\bf 2.)~Local rescaling:}
The crucial case is when $x_i=(s_i,t_i)$ converges to $\pd\H^2$ so fast
that even $t_i R_i C_i \to 0$. So for simplicity we assume here that $x_i\in\pd\H^2$.
Then we can restrict $\X_i$ to half balls of radius $\d_i:=(C_i R_i)^{-1}\leq \ep_i$ 
and rescale them to connections $\tX_i(y):=\X_i((x_i,z_i)+\d_i y)$ 
on the half ball $D^4\subset\H^4$ of radius $1$ centered at $0$.
The rescaled connections then satisfy
\begin{equation} \label{nonzero}
|F_{\tA_i}(0)| \geq 1 .
\end{equation} 

\noindent
{\bf 3.) $\mathbf{L^p}$-decay of $\mathbf{F_\tX}$ for $\mathbf{p<3}$:}
By a calculation similar to lemma~\ref{lem L2 bound} 
for the curvature of the rescaled connections one obtains 
for all $2<p<3$
\begin{equation} \label{zero}
\|F_{\tX_i}\|_{L^p(D^4)} \to 0 .
\end{equation}

\noindent
{\bf 4.) $\mathbf{\cC^0}$-estimates for $\mathbf{F_\tA}$ in terms of 
$\mathbf{L^p}$-bounds on $\mathbf{F_\tX}$ for $\mathbf{p>\frac 83}$:}
From (\ref{zero}) for $p>2$ and Uhlenbeck's weak compactness theorem~\ref{thm:weak comp}
we know that (up to gauge and taking a subsequence) the rescaled connections
$\tX_i\in\cA(D^4)$ converge to a flat connection in the weak $W^{1,p}$-topology.
One obtains stronger estimates from the fact that the rescaling preserves 
the anti-self-duality equation.
This implies $\cC^\infty$-convergence of the $\tX_i$ away from the boundary $\pd\H^4$.
At the boundary, the local rescaling has lost the global part of the Lagrangian
boundary condition, but the slice-wise flatness persists,
$ F_{\tX_i}|_{\{(s,0)\}\times\R^2} = 0 $.
With this one can go through the steps b)--e) in section~\ref{sec:C}
to obtain $W^{2,q}$-estimates on some components of the $\tX_i$.
One then feeds these back into c) and d) to obtain
$W^{2,q}$-bounds on $\nabla\tilde\P_i$ and $\nabla\tilde\Psi_i$, where the derivative
$\nabla$ is only in the $\R^2$-directions corresponding to~$\rT\S$. 

Now we need to assume (\ref{zero}) with $p>\frac 83$, then we can work with $q>2$ and 
the above bounds are just strong enough to imply $\cC^0$-convergence of the curvature part 
$*F_\tA = \pd_t\tilde\P - \pd_s \tilde\Psi + [\tilde\Psi,\tilde\P]$.
Since this convergence is to a flat connection, it provides a contradiction to
(\ref{nonzero}).
Note that this contradiction between 3) and 4)
crucially relies on the celebrated fact 
\fbox{$\,\frac 83<3$}.

\subsection{Proof of Removability of Singularities} \hspace{2mm} 

The proof of theorem~\ref{thm:RemSing} in case (in-L,boundary)
proceeds through the subsequent three propositions.
Throughout we denote by $D_r:=D_r(0)\subset\H^2$ the half ball of radius $r>0$, 
by $D_r^*:=D_r\setminus\{0\}$ the punctured half ball, and we will use 
polar coordinates $(r,\p)\in D_1^*$ with $r\in(0,1]$ and $\p\in[0,\pi]$.

We will consider solutions of (\ref{bvp}) on $D_1^*\times\S$, 
that is anti-self-dual connections which satisfy the Lagrangian
boundary condition on $\{(s,0)\}\times\S$ for $s\neq 0$.
An important tool for a connection $\X\in\cA(D_1^*\times\S)$ with finite energy 
$\int_{D_1^*\times\S}|F_\X|^2<\infty$ is its
{\bf energy function} $\cE:(0,1]\to[0,\infty)$ given~by 
\begin{align*}
\cE(r):=\; \half \int_{D_r^*\times\S} |F_\X|^2  
\;\; \biggl[ \,=\; \lim_{\d\to 0} \half \int_{(D_r\setminus D_\d)\times\S} |F_\X|^2 
\;=\; \lim_{\d\to 0} \bigl( \cE(r) - \cE(\d) \bigr). \biggr]
\end{align*}
The above calculation shows that finite energy directly
implies ${\cE(\d)\to 0}$ as $\d\to 0$.
For a finite energy solution of (\ref{bvp}) one thus obtains mean value inequalities
as in section~\ref{mean} on sufficiently small punctured balls.

\begin{prp} \label{thm B} \cite[Lemma~5.4]{W bubbling}
There are constants $C$ and $\ep>0$ such that the following holds.
Let $\X\in\cA(D_1^*\times\S)$ be a solution of (\ref{bvp}) and suppose that
$\cE(2r)\leq\ep$ for some $r\in(0,\half]$. Then for all $\p\in[0,\pi]$
\begin{enumerate}
\item
$\|F_\X (r,\p)\|_{L^2(\S)} \leq C r^{-1}\sqrt{\cE(2r)}$ ,
\item
$\|F_\X (r,\p)\|_{L^\infty(\S)} \leq C (r\sin\p)^{-2}\sqrt{\cE(2r)}$ .
\end{enumerate}
\end{prp}

\noindent
{\bf Sketch of Proof:}
The estimate (ii) is the mean value inequality for ${e=|F_\X|^2}$ that follows from
proposition~\ref{prp:mean}. Since lemma~\ref{lem:Lap} does not provide a control
on $\frac\pd{\pd\n}e$ we can only work on balls that are entirely contained 
in $D_{2r}^*\times\S$. 
When centered at $(r,\p,z)\in D_1^*\times\S$, their maximal radius is $r\sin\p$.

Next, write the connection as $\X=\P\ds+\Psi\dt+A$.
For the curvature component $F_A$, which vanishes at the boundary $\p\in\{0,\pi\}$, 
we can improve (ii) to $\|F_A(r,\p)\|_{L^\infty(\S)}\leq C r^{-2}$.
This follows from $\|F_\X(r,\p)\|_{L^2(\S)}\leq C r^{-1}$
similar to lemma~\ref{crucial} 
($\|F_A\|_{L^\infty(\S)}$ is essentially bounded by $\|F_\X\|_{L^2(\S)}^2$). 
The latter estimate is proven by an indirect argument as in section~\ref{eq}.
This uses the mean value inequality for $e=\|F_\X\|_{L^2(\S)}^2$ from
proposition~\ref{prp:mean}, based on lemma~\ref{lem:Lap} and again lemma~\ref{crucial}.

Once $\|F_A(r,\p)\|_{L^\infty(\S)}\leq C r^{-2}$ is established that way, one can
use it again in the mean value inequality for $e=\|F_\X\|_{L^2(\S)}^2$. It provides
$\laplace e \leq C r^{-2} e$ on (partial) balls of radius $\half r$ around $(r,\p)$.
The claim (i) then follows directly.\\

The curvature decay established here is almost sufficient to remove the singularity. 
The exponent of $r$ only has to be slightly improved to achieve the conditions
in the following removable singularity result. 
This improvement will finally be achieved in the crucial proposition~\ref{thm A} 
by a control on the speed of convergence of the energy function $\cE(r)\to 0$ as $r\to 0$.

\begin{prp} \label{thm C} \cite[Theorem~5.3]{W bubbling}
Let $\X\in\cA(D_1^*\times\S)$ and suppose that for some
constants $C$ and $\b>0$ and for all $(r,\p)\in D_1^*$
\begin{enumerate}
\item
$\|F_\X (r,\p)\|_{L^2(\S)} \leq C r^{\b-1}$ ,
\item
$\|F_\X (r,\p)\|_{L^\infty(\S)} \leq C (\sin\p)^{-2} r^{\b-2}$ .
\end{enumerate}
Then there exists $p=p(\b)>2$ and 
a gauge transformation $u\in\cG^{2,p}(D_1^*\times\S)$ such that
$u^*\X$ extends to a connection $\tX\in\cA^{1,p}(D_1\times\S)$.

Moreover, if $\X$ is a solution of (\ref{bvp}), then $\tX$ automatically
solves (\ref{bvp}) on $D_1\times\S$. A further gauge tranformation 
then makes $\tX\in\cA(D_1\times\S)$ smooth.
\end{prp}

\noindent
{\bf Sketch of Proof:}
To control the connection in terms of its curvature we fix 
a special gauge: Trivializing the bundle along rays $0<r\leq 1$ for
fixed $\p=\frac\pi 2$ and $z\in\S$ and then along $0\leq\p\leq\pi$ for 
fixed $r$ and $z\in\S$ we obtain
$$
\X = A + R \,\rd r + 0 \,\rd\p \qquad\text{with}\; R|_{\p=\frac\pi 2}= 0 .
$$
Here $A:D_1^*\to\Om^1(\S;\su(2))$ and $R:D_1^*\to\Om^0(\S;\su(2))$.
In this gauge we have $|\pd_r \X|\bigr|_{\p=\frac\pi 2}\leq |F_\X|$
and $|\pd_\p \X|\leq r |F_\X|$ since the curvature decomposes as
$$
|F_\X|^2 = |F_A|^2 + |\pd_r A -\rd_A R|^2 + r^{-2} |\pd_\p R|^2 + r^{-2} |\pd_\p A|^2 .
$$
The bounds (i) and (ii) combine to $|F_\X|\in L^p(D_1\times\S)$ for
some $p>2$ that only depends on $\b$. Roughly, they also imply
$\X|_{\{(r,\frac\pi 2)\}\times\S} \to A_0 \in \cA^{0,p}(\S)$ 
and $\X|_{\{r\}\times[0,\pi]\times\S} \to A_0\in\cC^0([0,\pi],\cA^{0,p}(\S))$ 
as $r\to 0$, and $A_0$ provides the extension over $\{0\}\times\S$.
In practice one constructs a family of connections $(\X_\ep)_{\ep\geq 0}$ 
on $D_1\times\S$ that coincide with $\X$ outside of $D_{2\ep}\times\S$ and equal to 
$A(\ep,\frac\pi 2)$ on $D_\ep\times\S$. 
Using (i) and (ii) this cutoff construction can be done such that
$\|F_{\X_\ep}-F_\X\|_{L^p(D_1\times\S)}\to 0$ as $\ep\to 0$.

By Uhlenbeck's compactness theorem~\ref{thm:weak comp} one then finds a sequence $\ep_i\to 0$
and gauge transformations $u_i\in\cG(D_1\times\S)$ such that
$u_i^*\X_{\ep_i}$ converges $W^{1,p}$-weakly to a limit connection $\tX\in\cA^{1,p}(D_1\times\S)$.
Note that on every compact subset of $D_1^*\times\S$ the sequence $\X_{\ep_i}$ eventually 
coincides with $\X$. So the above convergence also implies that (for a subsequence) the gauge
transformations $u_i$ converge to a limit $u\in\cG^{2,p}_{\rm loc}(D_1^*\times\S)$
in the weak $W^{2,p}$-topology on every compact set.
Then by the uniqueness of the limit $u^*\X=\tX|_{D_1^*\times\S}$, so $\tX$ is the claimed
extension.
If moreover $\X$ and hence $\tX$ are solutions of (\ref{bvp}) then the regularity 
theorem~\cite[Theorem~A]{W elliptic} for this boundary value problem asserts that
$\tX$ is gauge equivalent to a smooth solution.

\begin{prp} \label{thm A} \cite[Lemma~4.1]{W bubbling}
Let $\X\in\cA(D_1^*\times\S)$ be a solution of (\ref{bvp}) with finite energy 
$\cE(1)<\infty$.
Then for all $r\in(0,1]$
$$
\cE(r) \leq r^{\frac 1\pi} \cE(1) .
$$
\end{prp}

\noindent
{\bf Sketch of Proof:}
By the anti-self-duality equation $*F_\X=-F_\X$ we have
\begin{align*}
\half \int_{(D_{r_0}\setminus D_\d)\times\S} \la F_\X \wedge * F_\X \ra 
= \half \int_{(D_{r_0}\setminus D_\d)\times\S} 
\rd \la \X \wedge ( F_\X - \tfrac 16 [\X\wedge\X] ) \ra .
\end{align*}
This converges to $\cE(r_0)$ as $\d\to 0$. On the other hand, Stokes' theorem expresses this
as integral over $\pd(D_{r_0}\setminus D_\d)\times\S$. Our goal is to rewrite it as 
$\cF(A_\d) - \cF(A_{r_0})$ for a functional $\cF$ depending on
${A_r:=A(r,\cdot):[0,\pi]\to\cA(\S)}$.
Here as in proposition~\ref{thm C} we work in the special gauge $\X=A+R\dr$.
Then proposition~\ref{thm B}~(i) gives 
$\|\pd_\p A_r\|_{L^2(\S)}\leq C\sqrt{\cE(2r)}\to 0$ as $r\to 0$,
so the paths $A_r$ are $L^2$-short paths connecting $A_r(0),A_r(\pi)\in\cL_H$.
These contribute to $\cF$ on the boundary components
$\{r_0\}\times[0,\pi]\times\S$ and $\{\d\}\times[0,\pi]\times\S$.
So it remains to deal with the boundary components\footnote{
One could eliminate these by gluing in paths $A'_r:[0,\pi]\to\cL_H$
in the Lagrangian connecting $A_r(0),A_r(\pi)\in\cL_H$.
This would reach the goal with a functional $\cF=\cF(A_r,A'_r)$.
For the subsequent argument however, the $L^2$-length of the path $A'_r$
has to be controlled by the $L^2$-distance of its endpoints.
The crucial point would be to establish this fact for paths in a fixed gauge orbit --
a subtle nonlinear $W^{1,2}$-estimate for gauge transformations.
}at $\p=0$ and $\p=\pi$.
We identify these with $([-r_0,-\d]\cup [\d,r_0])\times\S$ and glue in the domain
$([-r_0,-\d]\cup [\d,r_0])\times H$. 
Now extending the families $A_r(0),A_r(\pi)\in\cL_H$ by
$\tA_r(0),\tA_r(\pi)\in\cA_{\rm flat}(H)$ preserves the value of
$\int \la F_\X \wedge F_\X \ra$, and
\begin{equation} \label{energy}
\cE(r_0) = \CS(A_\d,\tA_\d) - \CS(A_{r_0},\tA_{r_0}).
\end{equation}
Here we introduce the Chern-Simons functional for a path ${A:[0,\pi]\to\cA(\S)}$ 
with $L^2$-close ends $\tA(0),\tA(\pi)\in\cL_H$
and extensions ${\tA(0),\tA(\pi)\in\Af(H)}$,
\begin{align*}
\CS(A,\tA) &= 
-\half \int_0^\pi\int_\S \la A \wedge \pd_\p A \ra 
+ \tfrac 1{12} \biggl[ \int_H \rd \la \tA(\p) \wedge [\tA(\p)\wedge\tA(\p)] \ra 
\biggr]_{\p=0}^{\p=\pi} \\
&= -\half \int_0^\pi\int_0^\p \int_\S \la \pd_\p A(\th) \wedge \pd_\p A(\p) \ra \,\rd\th \,\dph \\
&\quad - \tfrac 1{12} \int_H \la \bigl[ (\tA(0)-\tA(\pi)) \wedge (\tA(0)-\tA(\pi)) \bigr] 
                  \wedge \bigl(\tA(0)-\tA(\pi)\bigr) \ra  .
\end{align*}
This magic identity together with the special choice of extensions as in 
lemma~\ref{lem:ext}~(ii) allow us to obtain the isoperimetric inequality
\begin{align*}
\bigl| \CS(A_r,\tA_r) \bigr|
&\leq \half \left( \int_0^\pi \bigl\| \pd_\p A_r \bigr\|_{L^2(\S)} \,\dph \right)^2
 + \tfrac 1{12}  \left( \bigl\| \tA_r(0)-\tA_r(\pi) \bigr\|_{L^3(Y)} \right)^3  \\
&\leq \bigl( \half + \tfrac {C_H^3}{12} \bigl\| A_r(0)-A_r(\pi) \bigr\|_{L^2(\S)} \bigr)
 \left( \int_0^\pi \bigl\| \pd_\p A_r \bigr\|_{L^2(\S)} \,\dph \right)^2 .
\end{align*}
For sufficiently short $A_r$ this implies 
$|\CS(A_r,\tA_r)|\leq \pi \int_0^\pi \bigl\| \pd_\p A_r \bigr\|_{L^2(\S)}^2 $.
As seen before this converges to $0$ as $r=\d\to 0$, and moreover it is bounded
$\pi r \dot\cE(r)$. So (\ref{energy}) provides the differential inequality
$\cE(r)\leq \pi r \dot\cE(r)$.
Integrating $\frac{\rd}{\rd r}\ln\cE(r)\geq (\pi r)^{-1}$ 
then proves the claimed decay of $\cE(r)$.

 \bibliographystyle{alpha}

\end{document}